\documentclass[psamsfonts]{article}
 
 \advance\textwidth 1.3in
 \advance\textheight 1in
 \advance\topmargin-0.5in
 \advance\oddsidemargin-0.65 in
 \advance\evensidemargin-0.5 in

\usepackage{amsfonts}

 \title{ Classification of Simple  $C^*$-algebras of 
Tracial Topological Rank Zero   
\thanks{Research 
partially supported by NSF grants DMS 9801482.
         AMS 1991 Subject Classification Numbers:
                         Primary 46L05,
                                46L35.
                        Key words: Classification of Simple $C^*$-algebras,
						tracial topological rank zero.
                                \protect\\}}
\author{Huaxin Lin\\
 Department of Mathematics\\
University of Oregon\\
Eugene, Oregon 97403-1222}
\date{}
\begin{document}
\maketitle

\newcommand{\CA}{$C^*$-algebra}
\newcommand{\SCA}{$C^*$-subalgebra}
\newcommand{\aue}{approximate unitary equivalence}
\newcommand{\ayue}{approximately unitarily equivalent}
\newcommand{\mops}{mutually orthogonal projections}
\newcommand{\hm}{homomorphism}
\newcommand{\pisca}{purely infinite simple \CA}
\newcommand{\andeqn}{\,\,\,\,\,\, {\rm and} \,\,\,\,\,\,}
\newcommand{\QED}{\rule{1.5mm}{3mm}}
\newcommand{\morp}{contractive completely
positive linear map}
\newcommand{\asmorp}{asymptotic morphism}
\newcommand{\arrow}{\rightarrow}
\newcommand{\tdsum}{\widetilde{\oplus}}
\newcommand{\pa}{\|}  
\newcommand{\ep}{\varepsilon}
\newcommand{\id}{{\rm id}}
\newcommand{\aueeps}[1]{\stackrel{#1}{\sim}}
\newcommand{\aeps}[1]{\stackrel{#1}{\approx}}
\newcommand{\dt}{\delta}
\newcommand{\yu}{\fang}
\newcommand{\ca}{{\cal C}_1}
\newcommand{\Ad}{{\rm ad}}

\newtheorem{thm}{Theorem}[section]
\newtheorem{Lem}[thm]{Lemma}
\newtheorem{Prop}[thm]{Proposition}
\newtheorem{Def}[thm]{Definition}
\newtheorem{Cor}[thm]{Corollary}
\newtheorem{Ex}[thm]{Example}
\newtheorem{Pro}[thm]{Problem}
\newtheorem{Remark}[thm]{Remark}
\newtheorem{NN}[thm]{}
\renewcommand{\theequation}{e\,\arabic{section}.\arabic{equation}}

\newcommand{\Ik}{ {\cal I}^{(k)}}
\newcommand{\Iz}{{\cal I}^{(0)}}
\newcommand{\Ii}{{\cal I}^{(1)}}
\newcommand{\Ip}{{\cal I}^{(2)}}

\begin{abstract}
We give a classification theorem 
for unital separable simple nuclear \CA s with 
tracial topological rank zero which satisfy 
the Universal Coefficient Theorem.
We prove that if $A$ and $B$ are two such 
\CA s and 
$$
(K_0(A), K_0(A)_+, [1_A], K_1(A))
$$
$$
\cong (K_0(B), K_0(B)_+, [1_B], K_1(B)),
$$
then  $A\cong B.$
\end{abstract}

\section{Introduction}
There has been  rapid progress in recent years in the 
program of classification of nuclear \CA s, in particular,
in the case that \CA s are simple 
and have real rank zero.
For example,
Elliott and Gong (\cite{EG2}) show that
simple AH \CA s of real rank zero (with slow dimension growth)
can be classified (up to isomorphisms) by their scaled ordered $K_0$ and 
$K_1.$ Kirchberg and Phillips (\cite{K} and \cite{P1}) 
show that separable nuclear purely infinite simple 
\CA s satisfying the Universal Coefficient Theorem (UCT) can also be classified 
by their $K$-theory.
These two results are the highlights of the program of classification of nuclear 
\CA s initiated by G. A. Elliott (see \cite{Ell2}). More recent results
in the case of \CA s which have real rank one can be found 
in a paper of Elliott, Gong and Li (\cite{EGL}). 
Many other mathematicians have contributed to this program.
(See for example, \cite{BEEK}, \cite{BBEK}, \cite{D1}, \cite{D3}, \cite{DE2}, \cite{DL1}, \cite{DL2},
\cite{DG}, \cite{Ell1}, \cite{Ell2}, \cite{Ell3}, \cite{Ell4}, \cite{EG1}, 
\cite{EE}, \cite{EG2},  \cite{EGLP}, \cite{EGS}, \cite{ER}, \cite{ES}, \cite{G1}, \cite{G2}, 
\cite{JS}, \cite{BK1}, \cite{BK2}, \cite{Li1}, \cite{Li2},  \cite{Ln2}-\cite{LnI},
\cite{LnQ}-\cite{Lntor},
\cite{LP1}, \cite{LP2}, \cite{LP3},
\cite{LS}, \cite{LqP}, \cite{NT}, \cite{P1}, \cite{P2}, \cite{Ro1}, \cite{Ro2}, 
\cite{Ro3},
\cite{ER}, \cite{Su1}, \cite{Su2}, \cite{Th1}, \cite{Th2} - an incomplete list).
The application of the classification program began to spread
to different fields, in particular, to the study of dynamical systems. 
In this paper we will only consider the case that 
\CA s are simple and have real rank zero and stable rank one.

During the development of classification 
of nuclear separable purely infinite simple 
\CA s, a special class of inductive limits of finite direct sums of even Cuntz-algebras 
was first classified by R\o rdam (\cite{Ro1}). R\o rdam's paper 
kicked off a series of  work on classification of nuclear 
purely infinite simple \CA s. More  
general examples of purely infinite simple \CA s were classified, such as
simple \CA s which are inductive limits of $O_n$ tensored with $C(S^1)$
(\cite{LP1}) and 
those \CA s arising as crossed products of stable simple AF-algebra with a 
${\mathbb Z}$-action determined by an automorphisms (\cite{Ro2}; see also 
\cite{ER} and other papers listed above). These 
purely infinite simple \CA s exhaust all possible 
$K$-theory. Later Kirchberg and Phillips  classified 
separable nuclear purely infinite simple \CA s (satisfying the UCT)
without assuming an inductive limit structure or other 
special forms.

Naturally, it is important to classify a class of 
(unital) separable nuclear stably finite simple \CA s
without assuming they are inductive limits of certain simple forms.

In \cite{Lntaf} and \cite{Lntr}, we introduced the notion 
of tracial topological rank. Recall that a  simple unital \CA\, $A$ is said 
to have tracial topological rank zero (written $TR(A)=0$), if 
for any $\ep>0,$ any finite subset ${\cal F}\subset A$ and any
$a\in A_+\setminus\{0\},$ there exists a nonzero projection $p$
and a finite dimensional \SCA\, $B\subset A$ with $1_B=p$ such that

(1) $\|px-xp\|<\ep,$

(2) ${\rm dist}(pxp, F)<\ep$ for all $x\in {\cal F}$ and 

(3) $1-p$ is equivalent to a projection in ${\rm Her}(a)={\overline{aAa}}.$

Roughly speaking, if an AF-algebra could be described as \CA\, which can 
be approximated by
finite dimensional \SCA s in norm, then 
a \CA\, $A$ with $TR(A)=0$ could be described 
as \CA s which can be approximated by finite dimensional \SCA s
in ``measure," or rather in trace.
This also resembles a key structure of unital separable nuclear purely infinite
simple \CA\,: Every such \CA \, $A$ has the following property:
for any $\ep>0$ any finite subset ${\cal F}\subset A,$
there exists a nonzero projection $p$ and a \SCA\, $B\cong O_2$ such 
that

(1') $\|px-xp\|<\ep$ and 

(2') ${\rm dist}(pxp, B)<\ep$ for all $x\in {\cal F}.$

((3') Note that $1-p$ is equivalent to a projection in $Her(a)={\overline{aAa}}$ since
$A$ is purely infinite simple.) This property plays a very important  
role in the classification theorem of Kirchberg and Phillips.

It is shown (in \cite{Lntaf}) that a unital separable simple \CA\, $A$ with $TR(A)=0$
has real rank zero, stable rank one, weakly unperforated $K_0(A)$ 
and is quasidiagonal.
Simple AH-algebras with slow dimension growth have
stable rank one (see \cite{DNNP}) and are also quasidiagonal and have weakly unperforated
$K_0$-groups.
It is proved in \cite{EG2} that 
all simple AH-algebras with slow dimension growth and with real rank zero
have tracial topological rank zero. 
It seems that the class of simple \CA s with tracial topological rank zero
is  the right substitute for class of quasidiagonal simple \CA s with real rank zero,
stable rank one and weakly unperforated $K_0.$ 
It was shown in \cite{LnI} that unital separable simple \CA s
which are inductive limits of type I \CA s with real rank zero, 
stable rank one, weakly unperforated $K_0$ and unique tracial states 
have tracial topological rank zero. 
N. C. Phillips (\cite{LqP}) shows certain simple crossed products arising from smooth 
minimal diffeomorphisms have tracial topological rank zero.
In studying non-commutative shift, recently, A. Kishimoto (\cite{Ki}) 
shows that the crossed product of the two sided infinite tensor product 
$\otimes_{\mathbb Z} A$ by the shifts has tracial topological rank zero.

Tracial topological rank zero (or TAF) for simple \CA s was introduced 
for the purpose of applying  the so-called uniqueness theorem 
originally established in \cite{Lnsu}. An immediate consequence 
of the uniqueness theorem and the notion of TAF is that 
a unital separable simple nuclear \CA s of tracial topological rank zero
with the same $K$-theory as that of UHF-algebra $Q$ with $K_0(Q)={\mathbb Q}$
is actually isomorphic to $Q$ (\cite{Lntaf}). 
Further attempts were made to classify 
all separable unital simple nuclear \CA s with tracial topological rank zero.
For example, it was shown that if $A$ is a separable unital simple nuclear \CA s with 
$RR(A)=0$ and $K_0(A)={\mathbb Q}\oplus {\rm ker} \tau,$ where $\tau$ is  
the unique trace (see \cite{Lncltaf} and \cite{LnQ}), then it is isomorphic 
to any other such \CA s. This was also obtained in \cite{DE2}
in the case that $K_0(A)={\mathbb Q}.$ It is also shown 
in \cite{Lncltaf} that 
$A$ is a unital separable simple nuclear \CA\, with $TR(A)=0$
such that $K_0(A)$ is a subring of ${\mathbb Q}$ (such as 
${\mathbb Z}[1/2]$) and $K_1(A)$ is torsion free then 
it is isomorphic to a unital $A{\mathbb T}$-algebra with 
the same $K$-theory.

In this paper by combining several previous results, 
we give a classification theorem 
for unital separable simple nuclear \CA s with tracial topological rank zero 
which satisfy the UCT. \CA s in this class are classified (up to 
isomorphism) by 
their scaled ordered $K_0$-groups and $K_1$-groups.

The proof of classification theorem (\ref{IITM}) as before
consists of a uniqueness theorem, an existence theorem
and an intertwining argument of Elliott. 
With the results in \cite{Lnsu}  and the introduction 
of the notion of tracial topological rank (or TAF), 
the desired uniqueness theorem (\ref{IIIL3}) was established
in \cite{Lncltaf}. However, the version of the 
(half) existence theorem (\ref{IIIexold}) in \cite{Lncltaf} (see also \cite{DE2}) does not 
provide a true existence theorem for the purpose 
of classification of simple nuclear \CA s with tracial topological 
rank zero. Extra conditions 
on the order structure of $K_0(A)$ were needed in \cite{Lncltaf},
\cite{DE2} and \cite{LnQ}. 
The main technical result in this paper 
is to recover the order information lost 
in \ref{IIIexold}. 

In Section 2 we show that a nuclear separable 
simple \CA s with tracial topological rank zero
has an approximate structure which is similar 
to a standard inductive construction.
Section 3 gives a special version 
of an interpolation property. In section 
4 we establish the required existence theorem, by applying 
section 2 and 3 together with the half existence 
theorem (\ref{IIIexold} from \cite{Lncltaf}) and a cutting lemma
from our previous result in \cite{Lntor}.
After we establish the required existence theorem, 
the  classification theorem can be established 
using the results already established in \cite{Lncltaf}.

{\bf Acknowledgements} This papar was written   while the author was visiting East China 
Normal University in Summer 2000 and visiting Mathematical Science Research Institute
at Berkeley. Research at MSRI is supported in part by NSF grant DMS-9701755. 
It is also partially supported by a grant from the NSF,  
and the Zhi-Jiang
Scholarship at the ECNU. 
The author benefited from conversation with Guihua Gong and N. C. Phillips.
He would also like to thank G. A. Elliott for inspiration 
and encouragement during the work of this and related research.

\section{An Approximation Construction}

Recall that a unital simple \CA\, $A$ is said to have 
tracial topological rank zero (written $TR(A)=0$),
if for any finite subset ${\cal F}\subset A$ any $\ep>0$
and any nonzero element $a\in A_+,$ there exists 
a finite dimensional \SCA\, $B\subset A$ with $p=1_B$ 
such that

(1) $\|px-xp\|<\ep,$

(2) ${\rm dist}(pxp, B)<\ep$ for all $x\in {\cal F}$ and 

(3) $1-p$ is equivalent to a projection in ${\overline{aAa}}.$

\vspace{0.2in}

\begin{NN}\label{IC}
{\rm 
We begin by reviewing a construction of \CA s of real rank zero.
This construction appeared in several places (\cite{BBEK}, 
\cite{Go} and \cite{D1}). We present a construction of  Dadarlat's version. 
Recall that a \CA\, $B$ is said to be residually finite 
dimensional (RFD), if $B$ has a separating family 
of finite dimensional irreducible representations.
If $B$ is separable RFD \CA, then 
there is a sequence of finite dimensional irreducible 
representations $\{\pi_n\}$ of $B$ which separates $A.$
Let $B$ be such a RFD \CA, $\{\pi_n\}$ be 
a separating sequence of finite dimensional irreducible 
representations and $\{x_n\}$ be a dense sequence of 
the unit ball in $B.$ Suppose that 
the rank of $\pi_n$ is $r(n).$ 

Let $A_1=B.$ Set $R(2)=1+k(1)r(1)$ and $A_2=M_{R(2)},$
Define a \hm\, $h_1: A_1\to A_2$
by 
$$
h_1(a)={\rm diag}(a,\pi_1(a),...,\pi_1(a))\,\,\,{\rm for}\,\,\, a\in A_1,
$$
where $\pi_1$ repeats $k(1)$ times.
Let $A_3=M_{R(3)}(A_2),$ where $R(3)=1+k(2,1)r(1)+k(2,2)r(2),$ and define 
$h_2: A_2\to A_3$ by
$$
h_2(a)={\rm diag}(a, {\tilde \pi}_1(a),....,{\tilde \pi_1}(a),
{\tilde \pi_2}(a),...,{\tilde \pi_2}(a))
$$
for $a\in A_2,$ where ${\tilde \pi_1}=\pi_1\otimes {\id}_{k(1)r(1)},$
${\tilde \pi_2}=\pi_2\otimes \id_{k(1)r(1)}$
and ${\tilde \pi}_1(a)$ repeats $k(2,1)$ times and 
${\tilde \pi_2}(a)$ repeats $k(2,2)$ times.
Continuing this fashion, we obtain 
a sequence of \hm s $h_n: A_n\to A_{n+1}$ and set $A=\lim_n(A_n, h_n).$
If $k(n,i)\to\infty$ as $n\to \infty$ for each $i,$
as in \cite{Lntaf}, $A$ is simple and $TR(A)=0.$ 
Let $\phi_{n,n+m}=h_{m+n}\circ \cdots h_n.$
Then
$$
\phi_{1, n}(a)=\pmatrix{ a &&&&\cr
                        & \Phi_1^{(n)}(a) &&&\cr
					    && \Phi_2^{(n)}(a) &&\cr
						&&& \ddots & \cr
						&&&& \Phi_{n}^{(n)}(a) \cr}
$$
for $a\in A_1,$ where $\Phi_1^{(n)}(a)$ is  $s(1,n)$ many copies of $\pi_1(a),$
$\Phi_2(a)$ is $s(2,n)$ many copies of $\pi_2(a),$...,
$\Phi_{n}$ is $s(n,n)$ copies of $\pi_n(a).$
From the construction, $\inf_n\{s(k,n)r(k)/\sum_{i=1}^ns(i,n)r(i)\}$ is positive
for each $k$ and $\sum_{i=1}^ns(i,n)r(i)/(1+\sum_{i=1}^ns(i,n)r(i))\to 1.$  

It is clear that this construction can be made 
even more general (for example, $A_k$ is a finite direct sum 
of matrix algebra over $B$ with different size)  but 
still have the similar properties.  It is clear while this construction 
is special, it is typical. But much more is true. 

In what follows we will show that every nuclear 
separable simple \CA\, $A$ with $TR(A)=0$ has a similar 
structure, or at least an approximated one. 
We need to replace \hm s by approximate multiplicative morphisms
and $\pi_k$ will be replaced by those whose ranges 
are contained in finite dimensional \CA s. 
}
\end{NN}
\vspace{0.2in}

\begin{Lem}{\rm (cf. 5.3 in \cite{Lntaf})}\label{ILHn}
 Let $A$ be a separable unital nuclear simple \CA\, with $RR(A)=0.$
Let $\{A_n\}$ be a sequence of nuclear RFD \CA s such that 
$A_n\subset A_{n+1}$ and $A$ is the closure of $\cup_n A_n.$
Fix  $\ep_n>0$ which decreases to zero,
a sequence of nonzero \hm s $h_n$ from $A_n$ to a finite dimensional \CA\,
$F_n$ 
and a finite subset ${\cal F}_n\subset A_n.$  
There exists a sequence of non-zero projections $\{e_n\}$  
and a sequence of  monomorphisms $h_n': F_n\to A$ with $h'_n(1_{F_n})=e_n$ satisfying the following:

(1) $\|e_nx-xe_n\|<\ep_n$ and

(2) $\|h_n'\circ h_n(x)-e_nae_n\|<\ep_n$ for all $x\in {\cal F}_n.$

\end{Lem}

{\it Proof:} Fix an integer $n.$
Let $\ep>0$ and a finite subset ${\cal F}_n
=\{x_1,x_2,...,x_k\}\subset A_n$ be
given. 
Let $I={\rm ker} h_n\subset A_n$ and $B$ be the hereditary $C^*$-subalgebra
of $A$ generated by $I.$ 
Let $C$ be the closure of $A_n+B.$ 
We will show that for any $c\in C,
bc, cb\in B$ for all $b\in B.$ It suffices to consider 
those elements $c\in A_n.$
Note that $B$ is the closure
of $IAI.$ Fix $b\in B.$ For any $\ep>0,$ there is a positive
element $e\in I$ such that
$$
\|eb-b\|<\ep.
$$
We also have $ce\in I\subset B$ and $ceb\in B$ for all $c\in A_n.$
Therefore
$$
\|cb-ceb\|<\ep.
$$
This implies that $cb\in B.$ Similarly, $bc\in B.$
Thus  $C$ is  a $C^*$-subalgebra
containing $B$ as a (closed) ideal. Since $C/B=A_n/B\cap I=A_n/I\cong F_n,$
we see that $C/B\cong F_n.$
To save notation, we may write $F_n=C/B.$ 
Let $\pi: C\to C/B$ be the quotient map.
Hence we may identify $h_n$ with $\pi|_{A_n}.$

By 5.2 in \cite{Lncltaf}, every projection in $F_n$ lifts to a projection in $C.$
Note that $B$ has real rank zero, so
it admits an approximate identity consisting of projections.
Now a standard argument (see Lemma 9.8 in \cite{Ef}) shows that
$C$ contains a finite dimensional $C^*$-subalgebra $F_n'$ such that
$F_n'\cong F_n$ and 
$\pi(F_n')=F_n.$  
Thus there is a monomorphism ${\tilde h_n}: F_n\to F_n'$
such that $\pi\circ {\tilde h_n}= \id_{F_n}.$ 
In particular, ${\tilde h_n}\circ \pi(x)-x\in B$ for all $x\in C.$  

Let $q=1_{F_n'},$ then $(1-q)C(1-q)=(1-q)B(1-q).$
Write $F_n'=M_{n_1}\oplus M_{n_2}\oplus\cdots \oplus M_{n_k}$ and assume
that $d_i\in M_{n_i}$ are minimal projections in $M_{n_i}.$
Let $\{e_0^{(m)}\}$ be an approximate identity for
$(1-q)B(1-q)$ consisting of projections. Let
$\{c_i^{(m)}\}$ be an approximate identity for
$d_iBd_i$ consisting of projections, $i=1,2,...,k.$
Write the hereditary $C^*$-subalgebra of $C$
generated by $M_{n_i}$ as $M_{n_i}(d_iCd_i).$
Put $e_i^{(m)}=diag(c_i^{(m)},c_i^{(m)},...,c_i^{(m)})$
($c_i^{(m)}$ repeats $n_i$ times) and $E_m=e_0^{(m)}+\sum_{i=1}^ke_i^{(m)}.$
Then it is clear that $\{E_m\}$ forms an approximate
identity for $B$ consisting of projections.
Furthermore, by construction,
$$
E_mx=xE_m
$$
for each $x\in F_n'.$
Since ${\tilde h_n}\circ h_n(x_i)-x_i\in B,$ 
for sufficiently large $m,$
$$
\|(1-E_m)({\tilde h_n}\circ h_n(x_i)-x_i)(1-E_m)\|<\ep/2,\,\,\,\, i=1,2,...,k.
$$
Moreover, since $E_m ({\tilde h_n}\circ h_n(x_i))=({\tilde h_n}\circ h_n(x_i))E_m,$  
$$
\|(1-E_m)x_i-x_i(1-E_m)\|<\ep/2,\,\,\,\, i=1,2,...,k,
$$
and  $h_n'(a)=(1-E_m){\tilde h_n}(a)(1-E_m)$ ($a\in F_n$) defines 
a \hm. 
Let $e_n=(1-E_m).$ Then 
$h_n'(1_{F_n})=e_n,$ 
$$
\|e_nx_i-x_ie_n\|<\ep/2\andeqn 
\|h_n'\circ h_n(x_i)-e_nx_ie_n\|<\ep
$$
$i=1,2,...,k.$
\QED

\begin{Def}\label{IDmul}
Let $A$ and $B$ be \CA s, let $L: A\to B$ be a \morp, let $\ep>0$ and let
${\cal F}\subset A$ be a subset. $L$ is said to be
${\cal F}$-$\ep$-multiplicative,
if
$$
\|L(xy)-L(x)L(y)\|<\ep
$$
for all $x,y\in {\cal F}.$
\end{Def}

\begin{NN}\label{INN}

{\rm Let $A$ be a separable unital nuclear simple \CA\, with $TR(A)=0.$
It follows from a result of B. Blackadar and E. Kirchberg 
(\cite{BK1} and \cite{BK2}) that $A$ is an inductive limit  
of RFD \CA s. 
Let $\{A_n\}$ be a sequence of nuclear RFD \CA s such that 
$A_n\subset A_{n+1}$ and $A$ is the closure of $\cup_n A_n.$
Fix a finite subset ${\cal F}_1\subset A_1,$  $\dt_0>0$ and a \hm\, $h$ from $A_1$ 
to a finite dimensional \SCA\, $F_0.$  Suppose that 
$\{x_1,...,x_n, ...\}$ is a dense sequence of elements in the unit ball of $A.$
For the convenience, we assume that $x_1\in {\cal F}_1.$ 
Applying Lemma \ref{IC},
we obtain a nonzero projection $e_0>0$ 
and a \hm\, $h':  F_0\to A$ with $e_0=h'(1_{F_0})$ such that

($i_0$) $\|e_0a-ae_0\|<\dt_0/8$   and 

($ii_0$) $\|h'\circ h(a)-e_0ae_0\|<\dt_0/32$ for all $a\in {\cal F}_1.$
Set $H=h'\circ h.$

Since $TR((1-e_0)A(1-e_0))=0,$ 
there is a  
projection $q_1'\le (1-e_0)$ and a finite dimensional \SCA\, $F_1'$ with 
$1_{F_1'}=q_1'$ such that

($i_0'$) $\|q_1'x-xq_1'\|<\dt_0/8,$

($ii_0'$) ${\rm dist}(q_1'xq_1', F_1')<\dt_0/32$ for all $x\in {\cal F}_1$ and 

($iii'_0$) $\tau(1-q_1')<1/4$ for all tracial states on $A.$

Set $q_1=q_1'+e_0$ and $F_1=F_1'\oplus h(F_0).$ 
Thus

(i) $\|q_1x-xq_1\|<\dt_0/4,$

(ii) ${\rm dist}(q_1xq_1, F_1)<\dt_0/16$ for all $x\in {\cal F}_1$ and 

(iii) $\tau(1-q_1)<1/4$ for all tracial states $\tau$ on $A.$ 

Let ${\cal F}_2$ be the union of $\{x_2\}\cup {\cal F}_1$ and a set $S_1$ of standard generators
of $F_1.$ 
Since $A$ is nuclear, there is a \morp\, $L_1': q_1Aq_1\to F_1$ such that
$L_1'|_{F_1}={\id}_{F_1}.$ Set $L_1(a)=L_1'(q_1aq_1)$ for 
all $a\in A.$ Then $L_1$ is a \morp\, from $A$ to $F_1$ and $L_1|_{F_1}=\id_{F_1}.$  
Note that $L_1$ is $\{x_1\}$-$\dt_0/2$-multiplicative. 
There is $\dt_2'>0$ such that  for any ${\cal S}_1$-$\dt_2'$-multiplicative 
\morp\, $L$ from $F_1$ there is a \hm\, $h$ from $F_1$ such
that
$$
\|L(a)-h(a)\|<\dt_0/16 
$$
for all $a\in \{L_1(x_1)\}\cup S_1.$
Let $\dt_2=\min\{\dt_2', \dt_0/4\}.$ Since $TR(A)=0,$ there is a projection $q_2$ and 
a finite dimensional \SCA\, $F_2$ with $1_{F_2}=q_2$ such that

(i') $\|q_2x-xq_2\|<\dt_2/4$ for all $x\in {\cal F}_2,$ 

(ii') ${\rm dist}(q_2xq_2, F_2)<\dt_2/16$ for all $x\in {\cal F}_2$ and 

(iii') $\tau(1-q_2)<1/8$ for tracial states on $A.$

Since $A$ is nuclear, there is a \morp\, $L_2': q_2Aq_2\to F_2$ such 
that $L_2'|_{F_2}={\rm id}_{F_2}.$ Set $L_2(a)=L_2'(q_2aq_2)$ for all $a\in A.$
Note that $L_2$ is a \morp\, from $A$ to $F_2$ such that
$L_2|_{F_2}=\id_{F_2}$ and it is ${\cal F}_2$-$\dt_2$-multiplicative. 
Thus there is a \hm\, $h_2: F_1\to F_2$ such that 
$$
\|L_2(a)-h_2(a)\|<\dt_0/16
$$
for all $a\in \{L_1(x_1)\}\cup S_1.$  Since $L_2(z)\not=0$ 
(by (i'),(ii'))and (iii')) if $z\in S_1,$
$h_2$ has to be injective.

Continuing this way, we obtain a sequence of finite subsets 
${\cal F}_0, {\cal F}_1,...,{\cal F}_n,...$
in the unit ball of $A$ which is dense, a sequence of decreasing positive 
numbers $\dt_n$ (with $\dt_n<\dt_{n-1}/4^n$), a
sequence of projections $\{q_n\} \subset A,$ 
a sequence 
of finite dimensional \SCA s $F_n$ with $1_{F_n}=q_n$
and a sequence of \hm s $h_{n+1}\to F_n\to F_{n+1}$
such that

(1) $\|q_nx-xq_n\|<\dt_n/4$ for ${\cal F}_n$

(2) ${\rm dist}(q_nx_iq_n, F_n)<\dt_n/16,$ $i=1,...,n.$

(3) $\tau(1-q_n)<1/2^n$ for all tracial states $\tau$ on $A,$

(4) ${\cal F}_{n+1}\supset S_n,$ where $S_n$ is a set of standard generators of $F_n,$ 

(5) 
$
\|L_{n+1}(a)-h_{n+1}(a)\|<1/16^{n}
$
for all $a\in \{L_n({\cal F}_n)\}\cup\{S_n\},$
where $L_n:A\to F_n$ is a \morp\, such that 
$L_{n}|_{F_n}={\id}_{F_n},$ $n=1,2,....$

\vspace{0.1in}

Let $D_{n,1}, D_{n,2},..., D_{n,m(n)}$ be  simple summands of $F_n$ and 
$\pi_{n,i}: F_n\to D_{n,i}$
be the quotient map. Let $\Psi_n: A\to (1-q_n)A(1-q_n)$ defined by $\Psi_n(a)=(1-q_n)a(1-q_n)$ 
for $a\in A.$ Note that $\Psi_n$ are ${\cal F}_n$-$\dt_n/2$-multiplicative.
Define $J_n: A\to A$ by $J_n(a)=L_n(a)\oplus \Psi_n(a)$ (for $a\in A$).
Note that $J_n$ is ${\cal F}_n$-$\dt_n/2$-multiplicative.
Set $J_{m,n}=J_n\circ J_{n-1}\circ\cdots \circ J_m$
and $h_{m,n}=h_n\circ h_{n-1}\circ\cdots h_{m+1}:F_m\to F_n.$ Note also that
$J_{m,n}$ is ${\cal F}_m$-$\dt_m$-multiplicative.
To save notation, we will also use $L_n, \Psi_n,$ $J_n,$  $ J_{m,n},$ $h_m,$ and $h_{m,n}$ for 
$L_n\otimes {\id}_{M_k},$ $\Psi_n\otimes {\id}_{M_k}, J_n\otimes {\id}_{M_k}$
and $J_{m,n}\otimes {\id}_{M_k},$ $h_m\otimes {\id}_{M_k}$ and $h_{m,n}\otimes {\id}_{M_k}.$
}
\end{NN}

The following has been discussed in many recent papers (cf. \cite{Lncltaf}, \cite{Lntor}
and \cite{DE2}). We present here for the sake of notation. 
\vspace{0.2in}

\begin{Def}\label{KKL}
{\rm 
Let $A$ be a nuclear \CA\, and $B$ be a $\sigma$-unital \CA. Let 
${\cal T}(A,B)$ be the set of those extensions 
of $B$ by $A$ whose six-term exact sequence in $K$-theory  gives 
pure extensions 
$$
0\to K_i(B)\to K_i(E)\to K_i(A)\to 0\,\,\,\,\,\,\,i=0,1
$$
(i.e., every finitely generated subgroup of $K_i(B)$ splits).
Following R\o rdam, set\\
$KL(A,B)=Ext(A, SB)/{\cal T}(A, SB).$
Let $\beta\in KL(A,B),$ then $\beta$ gives an \hm\, from 
$K_0(A)\to K_0(B).$ Suppose that both $A$ and $B$ are stably finite.
Denote by $KL(A,B)_+$ those $\beta$ which induce  positive
\hm s.

Let $C_n$ be a commutative \CA\, with $K_0(C_n)={\mathbb  Z}/n{\mathbb Z}$
and $K_1(C_n)=0.$ Suppose that $A$ is a \CA.
Then $K_i(A, {\mathbb  Z}/k{\mathbb  Z})=K_i(A\otimes C_k).$
Let ${\bf P}(A)$ be the set of all projections in
$M_{\infty}(A),$ $M_{\infty}(C(S^1)\otimes A),$
$M_{\infty}((A\otimes C_m{\tilde)})$ and
$M_{\infty}((C(S^1)\otimes A\otimes C_m{\tilde )}).$
We have the following commutative diagram (\cite{Sc2}):
$$
\begin{array}{ccccc}
 K_0(A) &\to &K_0(A, {\mathbb  Z}/k{\mathbb  Z})
&\to& K_1(A)\\
\uparrow_{\bf k} & & & &\downarrow_{\bf k}\\
K_0(A) & \leftarrow & K_1(A, {\mathbb  Z}/k{\mathbb  Z})
& \leftarrow & K_1(A)
 \end{array}
$$
As in \cite{DL1}, we use the notation
$$
{\underline K}(A)=
\oplus_{i=0,1, n\in {\mathbb  Z}_+} K_i(A;{\mathbb Z}/n{\mathbb Z}).
$$
 By
$Hom_{\Lambda}({\underline K}(A),{\underline K}(B))$
we mean all \hm s from ${\underline K}(A)$ to ${\underline K}(B)$
which respect the direct sum
decomposition and the so-called Bockstein operations (see \cite{DL1}).
It follows from \cite{DL1} that if $A$ satisfies the Universal
Coefficient Theorem, then $Hom_{\Lambda}({\underline K}(A),{\underline K}(B))
=KL(A,B).$

Let $A$ and $B$ be two \CA s and $L: A\to B$ a completely positive
linear map. Then $L$ induces maps from $A\otimes C_m\to B\otimes C_m,$
from $C(S^1)\otimes A\otimes C_m$ to $C(S^1)\otimes B\otimes C_m,$
namely, $L\otimes {\rm id}.$ For convenience, we will also
denote the induced map by $L$.
Given a projection $p\in {\bf P}(A),$
if $L$ is ${\cal G}$-$\ep$-multiplicative
with sufficiently large ${\cal G}$ and sufficiently small
$\ep,$
$L(p)$ is close to a projection. Let $L(p)'$ be that projection.
Fix  finite subsets of ${\cal P}_1\subset {\bf P}(A).$
It is easy to see
that  $L(p)'$ and $L(q)'$ are in the same equivalence
class of projections of  ${\bf P}(A),$ if $p$ and $q$ are
in ${\cal P}_1$ and are
in the same equivalence class of projections
of ${\bf P}(A),$ provided that ${\cal F}$ is sufficiently large and
$\ep$ is sufficiently small. We use $[L](p)$ for the class of
projections containing $[L](p)'.$
In what follows, whenever we write
$[L](p),$ we assume that
${\cal F}$ is sufficiently large and $\ep$ is sufficiently small
so that $[L](p)$  is  well-defined
on ${\cal P}_1.$
Furthermore, abusing the language,
we write $[L]([p])$ as well as $[L](p),$ where
$[p]$ is the equivalence class containing $p.$
Suppose that $q$ is in ${\bf P}(A)$ with
$[q]=k[p]$ for some integer $k,$ by adding sufficiently many
elements (partial isometries) in ${\cal F},$ we can
assume that $[L](q)=k[L](p).$
Suppose that $G$ is a finitely generated group generated
by ${\cal P}$ and $G={\mathbb  Z}^n\oplus {\mathbb  Z}/k_1{\mathbb  Z}\oplus
\cdots {\mathbb  Z}/k_m{\mathbb Z}.$
Let $g_1,g_2,...,g_n$ be free generators of ${\mathbb Z}^n$
and $t_i\in {\mathbb  Z}/k_i{\mathbb Z}$ be the generator with order $k_i,$
$i=1,2,...,m.$ Since every element in $K_0(C)$ (for any unital \CA\, $C$) may be
 written as $[p_1]-[p_2]$ for projections $p_1, p_2\in A\otimes M_l,$
for some $l>0,$
with sufficiently large ${\cal F}$ and
sufficiently small $\ep,$ one can define
$[L](g_j)$ and $[L](t_i).$
Moreover (with sufficiently large ${\cal F}$ and
sufficiently small $\ep$), the order of $[L](t_i)$ divides $k_i.$
Then we can define a map $[L]|_{G}$
by
defining $[L](\sum_i^n n_ig_i+\sum_j^mm_jt_j)=
\sum_i^k n_i[L](g_i)+\sum_j^mm_j[L](t_j).$
Thus $[L]$ is a group \hm\, on $G.$
Note, in general, $[L]|_{{\cal P}}$ may {\it not}
coincide with $[L]|_{G}$ on ${\cal P}.$
However, if ${\cal F}$ is large enough and $\ep$ is
small enough, they coincide.
In what follows,
if ${\cal P}$ is given, we say $[L]|_{G}$ is well-defined  and write $[L]|_{G}$
if $[L]|_{\cal P}$ is well-defined, $[L]|_{G}$ is well-defined
and is a \hm\, and $[L]|_{{\cal P}}=[L]|_{G}$ on ${\cal P}.$
We will also use $[L]|_{\cal P}$ for $[L]|_{[{\cal P}]}$ whenever 
it is convenient.
}
\end{Def}

\begin{Def}\label{ITR}
  {\rm 
  Let $S$ be a compact convex set and ${\rm Aff}(S)$ be
 the space of all affine continuous functions on $S.$
${\rm Aff}(S)$ has the following order:
$$
{\rm Aff}(S)_+=\{f\in {\rm Aff}(S): f(s)>0\,\,\, s\in S\}\cup\{0\}.
$$
This is the order we use in this paper.

Let $A$ be a simple \CA\, with $TR(A)=0$ and with 
the tracial space $T(A).$ 
If $\tau\in T(A),$ we extend it to a trace on $A\otimes M_n$ 
on every $n$ by using $\tau\otimes Tr,$ where 
$Tr$ is the standard trace on $M_n.$ 
We will use $\tau$ for the extension.

Denote by $\rho_A: K_0(A)\to {\rm Aff}(T(A))$ the 
\hm\, defined by $\rho_B([p])(\tau)=\tau(p)$ 
for all $\tau\in T(A)$ and  for projections $p\in A\otimes M_n$
($n=1,2,....$).
We will often write $\tau([p])$ for $\rho_A([p]).$
}
\end{Def}

We now return to the construction in \ref{INN}.
\vspace{0.2in}

\begin{Lem}\label{ILT}
Let ${\cal P}\subset M_k(A)$ be a finite subset of projections. Assume 
that ${\cal P}\subset M_k(A_1)$ and ${\cal F}_1$ is sufficiently large and 
$\dt_0$ is sufficiently small.
 With notation in \ref{INN}, we may assume that 
$[L_{1+n}\circ J_{1,n}]|_{{\cal P}}$ and $[L_{n+1}\circ J_{1,n}]|_{G}$ are well defined,
where $G$ is the subgroup determined by ${\cal P},$ and 
$$
\lim_{n\to\infty}\tau([L_{1+n}\circ J_{1,n}]([p]))=\tau([p])
$$
and the convergent is uniform on $T(A).$ Furthermore,
$$
|\tau(h_{1+n}\circ\cdots\circ h_1([p]))-\tau(h_n\circ \cdots \circ h_1([p])|<1/2^{n+1}
$$
for all $\tau\in T(A)$ and 
$$
  \lim_{n\to\infty}\tau(h_{1+n}\circ\cdots h_{1}([p])\ge 
  (1-\sum_{k=1}^{\infty}2^{k+1})\tau(h_{1}([p]))>0.
$$
for $p\in {\cal P}$ and for all tracial states on $A$ 
(and the convergent is uniform on $T(A)$).
\end{Lem}

{\it Proof:}
Choose a finite subset ${\cal F}_1\subset A_1$ and $\dt_0>0$ so that $[L]|_{\cal P}$ and 
$[L]|_{G}$ are well defined provided 
that $L$ is a  ${\cal F}_1$-$\dt_0$-multiplicative \morp.
Without loss of generality, we may assume that  
$\dt_0<1/2.$ 
This implies that $[p]=
[J_{1,n}(p)]$ for $p\in {\cal P}.$ In particular 
$\tau([p])=\tau([J_{1,n}(p)])$ for  $p\in {\cal P}.$ 
The first limit formula together with the uniformness of the convergence 
follows from the fact that 
$\tau(1-q_n)<1/2^{n+1}$ 
for all tracial states $\tau$ 
and Theorem 6.8 in \cite{Lntr} (see also 3.4 and 3.6 in \cite{Lntaf})
that $A$ has the fundamental comparison property of Blackadar.
Second formula also follows from the the same fact that $\tau(1-q_k)<1/2^{k+1},$ $k=1,2,....$ 
\QED

\begin{NN}\label{IPoint}
{\rm 
Fix a finite subset ${\cal P}$ of projections of $A.$ 
Let $N$ be a sufficiently large integer.
Then 
\begin{eqnarray*}
[L_{N+1}\circ J_{N}]&=&[L_{N+1}\circ L_N]\oplus [L_{N+1}\circ \Psi_N]\\
&=&[h_{N+1}\circ ({\rm diag}(\pi_{N,1}, 
\pi_{N,2},...,\pi_{N,m(N)}))\circ [L_N]\oplus [L_{N+1}\circ \Psi_N]
\end{eqnarray*}
on ${\cal P},$ and 
 \begin{eqnarray*}
[L_{N+2}\circ J_{N, N+1}]
&=&[L_{N+2}\circ L_{N+1}\circ J_N]\oplus [L_{N+2}]\circ [\Psi_{N+1}\circ J_N]\\
&=&[L_{N+2}\circ L_{N+1}\circ L_N]\oplus [L_{N+2}\circ L_{N+1}\circ \Psi_N]
\oplus [L_{N+2}]\circ [\Psi_{N+1}\circ J_N]\\
&=&[h_{N,N+2}\circ({\rm diag}(\pi_{N,1},  \pi_{N,2},...,\pi_{N,m(N)}))\circ [L_N]\\
&&\hspace{0.3in} \oplus [L_{N+2}\circ L_{N+1}\circ \Psi_N]\oplus [L_{N+2}]\circ 
[\Psi_{N+1}\circ J_N].
\end{eqnarray*}
Moreover, 
 \begin{eqnarray*}
[L_{N+n}\circ J_{N, N+n-1}]
&=&[h_{N,N+n}\circ({\rm diag}(\pi_{N,1},  \pi_{N,2},...,\pi_{N,m(N)}))\circ [L_N]\\
&& \oplus [L_{N+n}\circ \Psi_{N+n-1}\circ J_{N, N+n-2}]
\oplus
[L_{N+n}\circ L_{N+n-1}\circ \Psi_{N+n-2}\circ J_{N, N+n-3}]\\
&& \hspace{0.3in}\oplus\cdots \oplus [L_{N+n}\circ \cdots L_{N+1}\circ \Psi_N].
\end{eqnarray*}
on ${\cal P}.$  
Denote by ${\tilde \psi_{(N,i)}}=\pi_{N,i}\circ [L_N],$ 
${\tilde \psi_{(N+1,i)}}=\pi_{N+1,i}\circ L_{N+1}\circ \Psi_{N},
..., {\tilde \psi_{(N+n,i)}}=\pi_{N+n-1,i}\circ L_{N+n-1}\circ\Psi_{N+n-2}.$ 
Let $c_{(N+n,i,m)}(\tau)=\tau(h_{N+n, N+n+m}\circ {\tilde \psi_{N,i}}(1_A))>0.$ 
Moreover, we have computed 
(in \ref{ILT}) that $c_{(N+n,i,m)}(\tau)\to c_{(N+n,i)}(\tau) >0$ uniformly (as $m\to\infty$)
on $T(A).$ 

Rearranging $\{{\tilde \psi_{N+m,i}}\}$ ($m=1,2,....$) as $\psi_j,$ $j=1,2,...$
Set $s(k)=\sum_{l=1}^km(N+l-1).$ Suppose that the image of $\psi_j$ is isomorphic to $M_{r(j)}.$
So the rank of the image of $\psi_j$ is $r(j).$ 

}
\end{NN}
 \begin{NN}\label{INP}
{\rm Identify $K_0(M_{r(j)})$ with ${\mathbb Z}.$  Let $g_j$ be the element in $K_0(A)$
 corresponding to the minimal projection in $\psi_j(A).$
 Fix $p\in {\cal P}.$ Set $z=tr\circ [\psi_j]([p]),$ where 
 $tr$ is the normalized trace on  matrix algebras.
Then $(r(j)z)g_j\in K_0(A).$ If $\psi_j={\tilde \psi_{N+m,i}},$ define
 $g_j^{(n)}=[h_{N+m, N+m+n}](g_j)$ and $a_j^{(n)}=c_{(N+m,i,n)}.$ Then
$\tau(g_j^{(n)})=a_j^{(n)}(\tau))/r(j).$
 }
 \end{NN}

We have the following lemma

\begin{Lem}\label{ILarr}
For any $p\in {\cal P},$
$$
\tau([p])=\lim_{k\to\infty}\sum_{j=1}^{s(k)}a_j^{(k)} tr\circ [ \psi_j]([p]) \,\,\,\,
{\rm uniformly\,\,\, on}\,\,\, T(A),
$$
where $tr$ denotes the normalized standard trace on matrix algebras.
Furthermore,  $a_j^{(k)}>0$ and $a_j^{(k)}\to a_j>0$ uniformly
on $T(A)$ (as $k\to\infty$),$j=1,2,...,s(k).$ 
\end{Lem}

{\it Proof:}
The expression of $\tau$ and the uniformness of the convergence
follow  from \ref{ILT} and the fact that
$\tau([\Psi_N+n](1))=\tau(1-q_{N+n})\to 0.$ The assertion that  
$a_j>0$ follows from the fact that each $c_{(m,i,n)}>0.$  
\QED

\begin{Cor}\label{ILlinf}
Let $G_0=G\cap K_0(A)$ and let ${\tilde \rho}: G\to l^{\infty}({\mathbb Q})$ be defined 
by
$$
g\mapsto (tr\circ [\psi_1](g),...,tr\circ[\psi_n],...).
$$
Then, ${\rm ker}{\tilde \rho}\subset {\rm ker}\tau.$
\end{Cor}

{\it Proof:} 
Clearly from \ref{ILarr}, if $tr\circ [\psi_j](g)=0$ for 
every $j,$ then $\tau(g)=0.$
\QED

\begin{NN}\label{IS}
{\bf Summery:}
{\rm The above construction (roughly) says, if $\ep>0,$ $\dt>0,$ 
finite subset of projections in $A_N\otimes M_k$ (for some $k>0$)
and finite subset  
${\cal F}\subset A_N$ are given, 
we may write ($a\in {\cal F}$) 
$$
a\approx_{\ep} \pmatrix{ \Psi_n(a) &&&&\cr
                         & \Phi_1^{(n)}(a) &&&\cr
						 & &\Phi_2^{(n)}(a) &&\cr
						 && & \ddots &\cr
                         &&&& \Phi_{s(n)}^{(n)}(a)\cr},
						 $$
where $\Phi_1^{(n)}=h_{N, N+n}\circ \psi_1,
\Phi_{s(k)}^{(n)}=h_{N+n-1, N+n}\circ \psi_{s(n)}.$
For any tracial state $\tau,$ we have 
$\tau(\Phi_j^{(n)})(a)=a_j^{n}(\tau)\cdot tr(\psi_1(a)).$
It is important that
$a_j^{(n)}$ converges to some positive element $a_j\in {\rm Aff}(T(A))$ 
uniformly on $T(A)$ and 
$$
|\tau([p])-\sum_{j=1}^{s(n)}a_j^{(n)} tr([\psi_j](p)|<\dt
$$
for all $\tau\in T(A)$ and for $p\in {\cal P}.$
}

\end{NN}

This concludes this section.



\section{A Technical Lemma}

\begin{NN}\label{IIAFF}
{\rm Let $S$ be a compact convex set.
Let $r$ be a positive integer or $r=\infty$ and 
${\mathbb D}\subset {\rm Aff}(S)$ be a subgroup. 
Set $l^{\infty}_r({\mathbb D})$  be set of 
$r$-tuple of elements in ${\mathbb D}$ (or the  set of bounded 
sequences in ${\mathbb D}$).

If $f=\{f_n\}\in l^{\infty}_r({\mathbb D}),$ we write
$$
\|f\|_{\infty}=\sup\{\|f_n\|: n=1,2,...,\}.
$$
}
\end{NN}

The following lemma will not be used in 
this paper. But it makes the next lemma much more transparent.

Let $\{x_{ij}\}$ be a $r\times \infty$ matrix, 
$a_j^{(n)}>0$ such that $a_j^{(n)}\to a_j>0,$ $j=1,2,....$ 
Suppose that 
$$
(x_{ij})_{r\times \infty} (a_j^{(n)})_{\infty\times 1}
\to z=(z_j)_{1\times \infty}.
$$
We want to find finitely many non-negative $b_1,...,b_n$ 
for some integer $n>0$ such that
$$
(x_{ij})_{r\times n}(b_j)_{n\times 1}=(z_j)_{n\times 1}.
$$
The following lemma states that much more than this is true.
Note that both constant $\dt$ and integer $K$ are important 
in the lemma.

\vspace{0.2in}

\vspace{0.2in}

\begin{Lem}\label{ILKEY}
Let $S$ be a compact convex set and ${\rm Aff}(S)$ be 
the space of all affine continuous functions on $S.$
Let ${\mathbb D}$ be a dense ordered subgroup of ${\rm Aff}(S).$
Let $\{x_{ij}\}_{0<i\le r,\,0<j<\infty}$ be a $r\times {\infty}$
matrix having  rank $r$ and with each $x_{ij}\in {\mathbb Q},$ and let 
$\{a_j^{(n)}\}$ be  sequences of positive  elements in ${\mathbb D}$
such that $a_j^{(n)}\to a_j(\, >0)$ uniformly on $S$ as $n\to\infty.$ 
For each $n,$
$$
(x_{ij})_{r\times n}v_n=y_n,
$$
where $v_n=(a_j^{(n)})_{n\times 1}$ is an $n\times 1$ column vector and 
$y_n=(b_i^{(n)})\in {\mathbb D}^r$ is an $r\times 1$ column vector.

Suppose that $y_n\to z$ for some $z=(z_j)_{r\times 1}\in {\mathbb D}^r$  
uniformly (in ${\rm Aff}(S)^r$ norm). 

Then there is $\dt>0$ and a positive integer $K>0$ satisfying the following:

For some sufficiently large $n,$
there is $u=(c_j)_{n\times 1}\in (1/K)({\mathbb D}^n_+)$ {\rm (}
where $c_j(\tau)>0$ for all $\tau\in S$ or
$c_j=0$ {\rm  )} such that
$$
(x_{ij})_{r\times n}u=z'
$$
if $z'\in {\mathbb D}^r_+$ and $\|z-z'\|_{\infty}<\dt.$ 
\end{Lem}

{\it Proof:}
To save notation, without loss of generality,
we may assume that
$
(x_{ij})_{r\times r}
$
has rank $r.$
Set $A_n=(x_{ij})_{r\times n}$ ($n\ge r$). 
Then there exists an invertible matrix 
$B\in M_r({\mathbb Q})$ (which does not depend on $n$)
such that
$BA_n=C_n,$ where $C_n=(c_{ij})_{r\times n},$  
$c_{ii}=1$ for $i=1,2,...,r,$ and $c_{ij}=0$ if $i\not=j,$ $j=1,2,...,r,$
and $c_{ij}\in {\mathbb Q}.$
Since $B\in M_r({\mathbb Q}),$ there is a positive integer $K>0$ such that
$K\cdot B\in M_r({\mathbb  Z}).$ 
Moreover, $K\cdot C_n\in M_{r\times n}({\mathbb Z}).$

Let $I_r$ be the $r\times r$ identity matrix.
We may write 
$$
C_n=(I_r, D_n'),
$$
where $D_n'$ is a $r\times (n-r)$ matrix.
Note that $K\cdot D_n'\in M_{r\times (n-r)}({\mathbb  Z}).$ 
Thus we have 
$$
C_nv_n=By_n
\andeqn I_rv_n'=By_n-D_nv_n,
$$
where $v_n'=(a_1^{(n)}, a_2^{(n)},...,a_r^{(n)})$ (as a column) and
$
D_n=(0, D_n')
$ 
is a $r\times n$ matrix.
Note that for any $n\times 1$ column vector $v$ with the form 
$(t_1, t_2,...,t_r,a_{r+1}^{(n)}, a_{r+2}^{(n)},...,a_n^{(n)}),$
$D_nv=D_nv_n.$
Since $a_j^{(n)}\to a_j>0$ uniformly on $S$ and 
$S$ is compact, there is $N_1>0$ such that 
$$
a_j^{(n)}(\tau)\ge \inf\{a_j(\tau)/2: \tau\in S\}>0 \,\,\,\,\,{\rm )}s\in S{\rm )}
$$
for all $n \ge N_1$ and $j=1,2,...,r.$
Let $0<\ep<min\{\inf\{a_j(\tau)\}/8: \tau\in S, \, j=1,2,...,r\}.$
There is $N_2>0$ such that 
$$
\|By_n-Bz\|_{\infty}<\ep/4
$$
if $n\ge N_2.$
There is $\dt>0$ depending only on $B$ 
such that, if $\|z-z'\|_{\infty}<\dt,$
$$
\|Bz-Bz'\|_{\infty}<\ep/4
$$
($B$ is determined by $\{x_{ij}\}$).
Therefore
$$
\|By_n-Bz'\|_{\infty}<\ep/2
$$
for all $n\ge N_2.$
Set $N=max\{N_1, N_2\}.$ 
Let 
$$
u'=Bz'-D_nv_n,
$$
where  $u'=(c_1,c_2,..,c_r)$ (column vector).
Since $I_rv_n'=v_n'$ and $I_ru'=u',$
we have
$$
\|u'-v'_n\|_{\infty}=\|(Bz'-D_nv_n)-(By_n-D_nv_n)\|<\ep
$$
if $n\ge N.$ 
Therefore $c_j>0$ for $j=1,2,...,r.$
Set $u=(c_1,...,c_r, a_{r+1}^{(n)}, a_{r+2}^{(n)},...,a_n^{(n)}).$
Then 
$$
I_ru'=Bz'-D_nu
$$
($n\ge N$). 
Since $I_ru'=Bz'-D_nu,$   $Bz'\in (1/K){\mathbb D}^r$ and 
$
K\cdot D_n\in M_{r\times(n-r)}({\mathbb Z})$
$$
(c_1,c_2,...,c_r)=u'=I_ru'\in (1/K)({\mathbb D}^r)_+.
$$ 
Since $D_n=C_n-I_r,$ we have 
$$
C_nu=Bz'.
$$
Finally, since $B$ is invertible, we have
$$
A_nu=z'.
$$

\QED 

\begin{NN}\label{Imatrix}
{\rm 
Let $G$ be a group and $A=(a_{ij})_{m\times k}\in M_{m\times k}({\mathbb Z}).$ 
Viewing  $y=(g_1,...,g_k)\in G^k$ as a column,
one defines $Ay= (\sum_{j=1}^ka_{1j}g_j, \sum_{j=1}^ka_{2j}g_j,...,\sum_{j=1}^ka_{mj}g_j).$
Thus $A$ maps $G^k$ to $G^m.$ 
If $A\in M_{m\times k}({\mathbb Z})$ and $B\in M_{r\times m}({\mathbb Z}),$
then $B(Ay)=(BA)(y).$ Note that if $B\in M_{r\times m}({\mathbb Q})$ but 
$BA\in M_{r\times k}({\mathbb Z})$
does not imply $B(Ay)=(BA)(y)$ in general, since $G$ might have torsion.  
}
\end{NN}

\begin{Lem}\label{ILKEY2}
Let $S$ be a compact convex set and ${\rm Aff}(S)$ be the set of 
all affine continuous functions on $S.$ 
Let ${\mathbb D}$ be a dense ordered subgroup of ${\rm Aff}(S)$ and
$G$ be an ordered group with the order determined 
by a surjective \hm\, $\rho: G\to {\mathbb D},$ i.e.,
$$
G_+=\{ g\in G: \rho(g)(\tau)>0\}\cup\{0\}.
$$
Let $\{x_{ij}\}_{0<i\le r,\,0<j<\infty}$ be a $r\times {\infty}$
matrix having  rank $r$ and with each $x_{ij}\in {\mathbb Q}_+,$ 
$r(j)\in {\mathbb N}$ such that $r(j)x_{ij}\in {\mathbb Z}_+$ for all $i,j,$ 
$g_j^{(n)}\in G$ such that
$\rho(g_j^{(n)})=a_j^{(n)}(\tau)/r(j)$ ($\tau\in S$), where 
 $\{a_j^{(n)}\}$ is a  sequence of positive  elements in ${\mathbb D}$
such that $a_j^{(n)}\to a_j(\, >0)$ uniformly on $S$ as $n\to\infty.$ 
For each $n,$ ($s(n)\ge n$)
$$
(r(j)x_{ij})_{r\times s(n)}{\tilde v}_n={\tilde y}_n,
$$
where ${\tilde v_n}=(g_j^{(n)})_{s(n)\times 1}$ is a $s(n)\times 1$ column  and 
${\tilde y_n}=({\tilde b_j^{(n)}})\in {G}^r$ is a $r\times 1$ column vector.
Set $b_i^{(n)}=\rho({\tilde b_j^{(n)}})$ and $y_n=(b_i^{(n)}).$ 
Suppose that $y_n\to z$ on $S$ uniformly on $S$
for some $z=(z_j)_{r\times 1}\in {\mathbb D}^r$  
(in ${\rm Aff}(S)^r$ norm). 

Then there is $\dt>0$ and a positive integer $K>0$ satisfying the following:

For some sufficiently large $n,$
there is $u=({\tilde c_j})_{s(n)\times 1}\in G^{s(n)}_+$ 
($\rho({\tilde c_j})>0$ or $\rho({\tilde c_j})=0$ ) such that
$$
(r(j)x_{ij})_{r\times s(n)}{\tilde u}={\tilde z}'
$$
if ${\tilde z}'\in G^r$  and there is ${\tilde z''}\in G^r$ such that
$K^3{\tilde z''}={\tilde z'}$ and $\|z-Mz'\|_{\infty}<\dt,$
where ${\tilde z}'=({\tilde z_1'},...,{\tilde z_r'})$ such that $z_j'=\rho({\tilde z}_j'),$
$j=1,2,...,r$ and $M$ is a positive integer. 
\end{Lem}

{\it Proof:} 
The proof is a repetition of the proof of \ref{ILKEY} with some  modification.
The important difference is that $G$ may have torsion.
We proceed as in the proof of \ref{ILKEY}. 
To save notation, without loss of generality,
we may assume that
$
(x_{ij})_{r\times r}
$
has rank $r.$
Set $A_n'=(r(j)x_{ij})_{r\times s(n)}\in M_{r\times s(n)}({\mathbb Z})$ ($s(n)\ge r$). 
Then there exists an invertible matrix 
$B'\in M_r({\mathbb Q})$ (which does not depend on $n$)
such that
$B'A_n'=C_n,$ where $C=(c_{ij})_{r\times s(n)},$  
$c_{ii}=1$ for $i=1,2,...,r,$ and $c_{ij}=0$ if $i\not=j,$ $j=1,2,...,r,$
and $c_{ij}\in {\mathbb Q}.$
Since $B'\in M_r({\mathbb Q}),$ there is a positive integer $K>0$ such that
$K\cdot B'\in M_r({\mathbb  Z}).$ 
Moreover, $K\cdot C_n\in M_{r\times s(n)}({\mathbb Z}).$ We may also assume that
$K\cdot (B')^{-1}\in M_r({\mathbb Z}).$ 

Let $I_r$ be the $r\times r$ identity matrix.
We may write 
$$
C_n=(I_r, D_n'),
$$
where $D_n'$ is a $r\times (s(n)-r)$ matrix.
Note that $K\cdot D_n'\in M_{r\times (s(n)-r)}({\mathbb  Z}).$ 
Since ${\mathbb D}$ is dense in ${\rm Aff}(S),$ there are $ \xi_n\in G^{s(n)}$
such that $\xi_n=({\tilde d}_j^{(n)})_{s(n)\times 1}$ and
$$
\|K^3\rho({\tilde d}_j^{(n)})-a_j^{(n)}/Mr(j)\|<{1\over{s(n)^2\cdot 2^n}},
\,\,\,\,j=1,2,...,s(n), n=1,2,....
$$
Let ${\tilde w_n}=K^3\xi_n$ and $\rho({\tilde d}_j^{(n)})=d_j^{(n)}.$
Then $K^3d_j^{(n)}\to K^3d_j=a_j/Mr(j)>0$  uniformly on $S.$  
Set ${\tilde y}_n'=A_n'{\tilde w}_n=K^3A_n'\xi_n.$ Note that $A_n'\xi_n\in G^r$
($A_n'\in M_{r\times s(n)}({\mathbb Z})$).
Let $\rho^{(s(n))}: G^{s(n)}\to {\mathbb D}^{s(n)}$ be defined 
by $\rho^{(s(n))}(g_1,...,g_{s(n)})=(\rho(g_1),...,\rho(g_{s(n)}))$ for $g\in G.$
Let $y_n'=\rho^{(s(n))}({\tilde y}_n').$ Then, from the construction above, 
$y_n'\to z/M$ uniformly on $S.$

We have 
$$
C_n{\tilde w_n}=K^3 C_n\xi_n=K^3B'A_n'\xi_n (=B'{\tilde y_n}')
\andeqn I_r{\tilde v_n'}=B'{\tilde y_n}'-D_n{\tilde w_n},
$$
where ${\tilde v_n'}=(K^3{\tilde d_1^{(n)}},...,K^3{\tilde d_r^{(n)}})$ 
 (as a column) and
$
D_n=(0, D_n')
$ 
is a $r\times s(n)$ matrix.
Note that for any $s(n)\times 1$ column vectors ${\tilde v}$  and $v$ with the form
$(s_1,...,s_r, K^3{\tilde d_{r+1}^{(n)}},K^3{\tilde d_{r+2}^{(n)}},..., K^3{\tilde d_{s(n)}^{(n)}})$ 
and\\ 
$(t_1, t_2,...,t_r,K^3d_{r+1}^{(n)}, K^3d_{r+2}^{(n)},...,K^3d_{s(n)}^{(n)}),$
$D_n {\tilde v}=D_n{\tilde w}_n$ and 
$D_nv=D_nw_n,$ where ${\rho^{(s(n))}}({\tilde v})=v$ and ${\rho^{(s(n))}}({\tilde w_n})
=w_n.$

Since $d_j^{(n)}\to d_j>0$ uniformly on $S,$  there is an $N_1>0$ such that 
$$
d_j^{(n)}\ge \inf\{d_j/2(\tau): \tau\in S\}>0
$$
for all $n \ge N_1$ and $j=1,2,...,r.$
Let $0<\ep<min\{\inf\{d_j(\tau): \tau\in S\}/8K^3: j=1,2,...,r\}.$
There is $N_2>0$ such that 
$$
\|B'y_n-B'z\|_{\infty}<\ep/2
$$
if $n\ge N_2.$
Thus there is $N_3>0$ such that
$$
\|B'(y_n')-B'(z/M)\|_{\infty}<\ep/2
$$
if $n\ge N_3$ ($y_n'\to z$ uniformly). 

There is $\dt>0$ depending only on $B'$ 
such that, if $\|z-Mz'\|_{\infty}<\dt,$
$$
\|B'y_n'-B'z'\|_{\infty}<\ep/2.
$$
if $n\ge N_3.$
Now let ${\tilde z}$ and ${\tilde z}''$ be as described in the lemma.
Set $N=max\{N_1, N_2, N_3\}.$ 
Let 
$$
B'{\tilde z'}-D_n{\tilde w_n}=K^3B{\tilde z''}-K^3D_n\xi_n=K^2u',
$$
where 
$u'=({\tilde c_1},{\tilde c_2},..,{\tilde c_r})\in G^r$ (column).
Set $u''=K^2u'.$ 
Let $\rho^{(r)}(u')=(c_1,c_2,...,c_r)\in {\mathbb D}^r$ We may write 
$$
I_ru''=B'{\tilde z'}-D_n{\tilde w_n}.
$$
Since $I_r{\tilde v_n'}={\tilde v_n'}$ and $I_ru'=u',$
we have
$$
\|\rho^{(r)}(u'')-\rho^{(r)}({\tilde v'_n})\|_{\infty}=
\|\rho^{(r)}(B'{\tilde z}'-D_n{\tilde w_n})-\rho^{(r)}(B'{\tilde y}_n'
-D_n{\tilde w_n})\|_{\infty}<\ep
$$
if $n\ge N.$ 
Since $\rho^{(r)}({\tilde v_n'})=(K^3d_1^{(n)}, K^3d_2^{(n)},...,K^3d_r^{(n)}),$ 
therefore $c_j(\tau)>0$ for all $\tau\in S$ and $j=1,2,...,r.$
Set 
$$
{\tilde u}=
(K^2{\tilde c_1},...,K^2{\tilde c_r}, K^3{\tilde d}_{r+1}^{(n)}, K^3{\tilde d}_{r+2}^{(n)},
...,K^3{\tilde d}_{s(n)}^{(n)})\andeqn
{\bar u}=({\tilde c_1},...,{\tilde c_r}, K{\tilde d}_{r+1}^{(n)},
...,K{\tilde d}_{s(n)}^{(n)}).
$$
Then ${\bar u}\in G^{s(n)}_+$ and
$$
I_ru''=K^3B'{\tilde z''}-K^2D_n{\bar u}=B'{\tilde z'}-D_n{\tilde u}
$$
($n\ge N$). 
Since $C_n=I_r+D_n,$ we have 
$$
C_n{\tilde u}=K^2C_n{\bar u}=K^2{\bar u}+ K^2D_n{\bar u}=
I_r{\tilde u}+D_n{\tilde u}
$$
$$
=I_ru''+D_n{\tilde w_n}=K^3B'{\tilde z''}=B'{\tilde z'}.
$$
We have $K (B')^{-1}, K B',K B'\in M_r({\mathbb Z}).$
Therefore 
$$
(B')^{-1}C_n {\tilde u}=K (B')^{-1} ( KC_n){\bar u}=K^2A_n{\bar u}=A_n {\tilde u}
$$
and
$$
(B')^{-1} B'{\tilde z}=(K\cdot (B')^{-1}) (K\cdot (B'))(K{\tilde z''})
=(K^2 I_r)(K{\tilde z''})={\tilde z'}.
$$
Hence 
$$
A_n{\tilde u}={\tilde z'}.
$$
\QED

\section{An Existence  Theorem}

The following theorem was proved 
in \cite{Lntaf}. A more general version 
was obtained by M. Dadarlat and S. Eilers (\cite{DE2}).

\begin{thm}{\rm (Theorem 5.9 in \cite{Lncltaf})}\label{IIIexold}
Let
$A$ be a separable \CA\, satisfying  UCT such that
$A$ is the closure of an increasing sequence $\{A_n\}$ of
RFD \CA s
and $B$ be  a unital nuclear separable
\CA. Then, for any $\alpha\in
Hom_{\Lambda}(\underline{K}(A),\underline{K}(B)),$ there exist
two sequences of completely positive contractions
$\phi_n^{(i)}: A\to B\otimes {\cal K}$ ($i=1,2$)
satisfying the following:

(1) $\|\phi_n^{(i)}(ab)-\phi_n^{(i)}(a)\phi_n^{(i)}(b)\|\to 0$ as
$n\to\infty,$

(2) for each $n,$ the images of $\phi_n^{(2)}$ are contained in
a finite dimensional $C^*$-subalgebra of $B\otimes {\cal K}$ and
for any finite subset ${\cal P}\subset {\bf P}(A),$
$[\phi_n^{(2)}]|_{\cal P}$ and
$[\phi_n^{(2)}]|_{G}$ are well defined for all
large $n,$
where $G$ is the subgroup generated by ${\cal P},$ 

(3) for  each finite subset of ${\cal P}\subset {\bf P}(A),$
there exists $m>0$ such that
$$
[\phi_n^{(1)}]|_{G}=\alpha+[\phi_n^{(2)}]|_{G}
$$
for all $n\ge m.$

(4) For each $n,$ we may assume that $\phi_n^{(2)}$ is
a \hm\, on $A_n.$

(The condition that $B$ is nuclear can be replaced by
the condition that each $A_n$ is nuclear).
\end{thm}

\vspace{0.2in}

It is shown in \cite{Lntaf} that if $A$ is a simple 
nuclear separable \CA s with $TR(A)=0,$ then 
$A$ is strong NF. 
It follows a result from B. Blackadar and E. Kirchberg (\cite{BK1} and \cite{BK2})
that every strong NF \CA\, is an inductive limit of residually finite dimensional
(RFD) \CA s. Thus the above theorem can be applied to 
simple nuclear separable \CA s with $TR(A)=0$ which satisfies the UCT. 

However the term $\phi_n^{(2)}$ prevents us to apply the theorem 
directly.

The following is proved in \cite{Lntor}:

\begin{Lem}{\rm (Lemma 1.8 in \cite{Lntor})}\label{ID} 
Let $A$ be a unital simple AH-algebra with 
slow dimension growth and with real rank zero.
Let $G_0$ be a finitely generated  subgroup
of $K_0(A)$ with decomposition
$G_0=G_{00}\oplus G_{01},$ where $G_{00}\subset {\rm ker}\rho_A$ and
$G_{01}$ is a finite generated free group such that
$(\rho_A)|_{G_{01}}$ is injective. Suppose that
${\cal P}\subset {\underline{K}}(A)$ is a finite subset
which generates a subgroup $G$ such that
$G\cap K_0(A)\supset G_0.$

Then, for any $\ep>0,$ any finite subset
${\cal F}\subset A,$
any $1>r>0,$ and any integer $K,$
there is an ${\cal F}$-$\ep$-multiplicative
map $L: A\to A$ satisfying the following:

(1) $[L]|_{\cal P}$ and $[L]|_{G}$ are  well-defined and $[L]|_{G}$
is positive on $G,$

(2) $[L]|_{G\cap {\rm ker}\rho_A}={\rm id}|_{G\cap {\rm ker}\rho_A},$
$[L]|_{G\cap K_0(A, {\mathbb Z}/k{\mathbb Z})}={\rm id}|
_{G\cap K_0(A, {\mathbb Z}/k{\mathbb Z})},$
$[L]|_{G\cap K_1(A)}={\rm id}|_{G\cap K_1(A)}$ and\\
$[L]|_{G\cap K_1(A, {\mathbb Z}/k{\mathbb Z})}=
{\rm id}|_{G\cap K_1(A, {\mathbb Z}/k{\mathbb Z})}$
for those $k$ with
$
G\cap K_i(A, {\mathbb Z}/k{\mathbb Z})\not=\emptyset
$ {\rm (}$i=0,1${\rm )},

(3) $\rho_A\circ [L](g)\le r\rho_A(g)$ for all $g\in G\cap K_0(A),$
 
(4) Let $g_1, g_2,..,g_l$ be positive generators
of $G_{01}.$
Then, there are $f_1,...,f_l\in K_0(A)_+$ such that
$$
g_i-[L](g_i)=Kf_i, \,\,\,i=1,2,...,l.
$$
\end{Lem}

Combining the above two results and \ref{ILKEY2}, we prove the 
existence theorem \ref{IILEX}.

To establish the following theorem, we first apply \ref{IIIexold} 
to obtain $\Phi_1.$ Then applying \ref{ID} to obtain $L$ 
with a small $r.$ 
We then try to construct $h$ with its range contained in 
a finite dimensional \SCA, and hope 
that $L\circ \Phi_1\oplus h$ will be the right map.
Given \ref{IIIexold}, the strategy is to construct $h$ which  
is a \hm\, on a finite dimensional \SCA.
So it is at least possible to have finite dimensional range and 
does not effect most other $KK$-theoretical information.
While this strategy sounds, the technical difficulties remain 
to be overcome. Immediately, one would need 
to extend some certain \hm\, from a subgroup of $K_0(C)$
($C$ is a finite dimensional \CA) to a \hm\, on $K_0(C).$
This surely requires some divisibility of the target group.
The integer $K$ in Lemma \ref{ID} comes to take care of that
problem.
However, the \hm\, on the group also has to be positive and 
the extension is also need to be positive.  
To make things worse, $\Phi_1$ does not even preserve
the order. Furthermore, one should also be careful when 
maps with finite dimensional ranges are used. They may not give 
trivial maps on (the part of) $K_0(A, {\mathbb Z}/k{\mathbb Z}).$
Lemma \ref{ILKEY2} together with a careful and extensive use 
of the integer $K$ in \ref{ID} will deal with these difficulties. 

Here is the existence theorem:

\vspace{0.2in}

\begin{thm}\label{IILEX}
Let $A$ be a unital separable simple nuclear \CA\, with $TR(A)=0$ 
which satisfies the UCT.
Suppose that $B$ is a unital simple AH-algebra with slow dimension growth, real rank zero
and 
$$
(K_0(B), K_0(B)_+, [1_B], K_1(B))\cong (K_1(A), K_1(A)_+, [1_A], K_1(A)).
$$
Let $\alpha\in KK(A,B)_+$ which carries the above isomorphism.
Then, for any finite subset ${\cal P}\subset {\bf P}(A),$ there is a sequence of 
\morp\, $H_n:A\to B$ such that

(i) $\|H_n(ab)-H_n(a)H_n(b)\|\to 0$ for all $a, b\in A$ as $n\to\infty$ and

(ii) $[H_n]|_{\cal P}=\alpha|_{\cal P}$ for all sufficiently large $n.$  
\end{thm}

{\it Proof:}
To save the notation, we may assume 
$$
(K_0(B), K_0(B)_+, [1_B], K_1(B))= (K_1(A), K_1(A)_+, [1_A], K_1(A)).
$$
So we will identify these two 4-tuples. 
Let ${\mathbb D}=\rho(K_0(A)).$ 
Note that ${\mathbb D}$ is a dense ordered subgroup of ${\rm Aff}(T(A)).$ 
 
Fix ${\cal P}$ and let ${\cal P}_0\subset {\cal P}$ be so that ${\cal P}_0$ generates 
a (finitely generated) subgroup $G_0$ so that $G_0=G({\cal P})\cap K_0(A),$ where 
$G({\cal P})$ is the subgroup generated by ${\cal P}.$
We may assume that ${\cal P}\subset {\bf P}(A_1).$ 
(Here we identify ${\cal P}$ with $j({\cal P}),$ where $j: A_1\to A$ is 
the embedding.) 
We may assume that $\{p_1,...,p_l\}={\cal P}_0,$ where $p_i\in M_k(A_1)$ are projections
(for some $k>0$).
Let ${\cal F}_0$ be a finite subset of $A_1$ and $\dt_0>0$ so that
any \morp\, $L$ from $A_1$ which is ${\cal F}_0$-$\dt_0$-multiplicative 
well defines $[L]|_{\cal P}$ and $[L]|_{G({\cal P})}.$  
Let $k_0$ be a positive integer such that
$$
G({\cal P})\cap K_i(A, {\mathbb Z}/k{\mathbb  Z})=\emptyset
$$
for all $k\ge k_0,$ $i=0,1.$

Step (I):  (Construct $\Phi_n$ and fix $H$)

It follows from \ref{IIIexold} that there is a sequence of \morp\,
$\Phi_n: A\to B\otimes {\cal K}$ such that 
$$
\|\Phi_n(ab)-\Phi_n(a)\Phi_n(b)\|\to 0
$$ for all $a, b\in A$ as $n\to\infty,$
 there is a sequence of \morp s $\{\Phi_n^{(0)}\}$ from $A$ to $B$ so that their  
images  contained in a finite dimensional \SCA s and
$\Phi_n^{(0)}|_{A_n}$ is a \hm, and 
$$
[\Phi_n]|_{G}=\alpha|_{G}+[\Phi_n^{(0)}]|_{G}\andeqn
[\Phi_n]|_{\cal P}=\alpha|_{\cal P}+[\Phi_n^{(0)}]|_{\cal P}
$$
for all $n.$ Fix $\Phi_1^{(0)}.$ Let $H: A_1\to A$ be a \hm\, such that
$H=h'\circ \Phi_1^{(0)},$ where $h'$ is a  monomorphism from the image of $\Phi_1^{(0)}$ to $A.$
The existence of such $h'$ is given by \ref{IC}.

For any finite subset ${\cal F}\subset A$ and $\eta>0$ (${\cal F}_0\subset {\cal F}$
and $\eta<\dt_0/2$), to save notation, we may 
assume that $\Phi_1$ is ${\cal F}$-$\eta$-multiplicative.
Fix this $H,$ ${\cal F}$ and $\eta/4$ (with $N=1$), 
let $q_n, \Psi_n, L_n, J_n, J_n, J_{m,n}, \psi_n, \dt_n$ (with 
$\sum_{n=1}^{\infty}\dt_n<\dt_0/2$) be as constructed in 
\ref{INN}. (The importance of $H$ will 
become clear later.)

Step (II): (Fix $M,$ $ K$ and $\dt>0$ )
Let ${\tilde \rho}: G({\cal P})\cap K_0(A)\to l^{\infty}({\mathbb Q})$ be defined by
$$
{\tilde \rho}([p_i])=(tr(\psi_1(p_i)),..., tr(\psi_m(p_i)),...)\,\,\,\,\,\,
{\rm for}\,\,\,
i=1,...,l.
$$ 
Note that ${\tilde \rho}(G_0)$ is a finitely generated 
torsion free group.  
Consider $\{{\tilde \rho}([p_1]), {\tilde \rho}([p_2]),...,{\tilde \rho}([p_l])\}.$
Suppose that the ${\mathbb Q}$-linear span of the above $l$ 
elements (in $l^{\infty}({\mathbb Q})$)
has rank $r.$ So a subset of 
$\{{\tilde \rho}([p_1]), {\tilde \rho}([p_2]),...,{\tilde \rho}([p_l])\}$
with $r$ elements has ${\mathbb Q}$- rank $r.$ 
Without loss of generality, we may assume that 
$ \{{\tilde \rho}([p_1]), {\tilde \rho}([p_2]),..., {\tilde \rho}([p_r])\}$ has rank 
$r$ and its 
${\mathbb Q}$-span 
includes all ${\tilde \rho}([p_i]),i=1,...,l.$
There is an integer $M>0$ such that for any $g\in {\tilde \rho}( G_0),$
$Mg$ is in the subgroup of ${\tilde G_0}$ generated by 
${\tilde \rho}([p_1]), {\tilde \rho}([p_2]),..., {\tilde \rho}([p_r]).$
Let $x_{ij}=tr({\tilde \psi_j}([p_i])),$ $i=1,2,...,r$ and $j=1,2,....$ So 
we may assume that 
$(x_{ij})_{r\times r}$ has ${\mathbb Q}$-rank $r.$ 
Let $g_j,$ $g_j^{(n)},$ $a_j^{(n)}$ and $a_j$ be as in \ref{INP}.
Let $z_i=\rho_A([p_i])\in {\mathbb D}$ and 
$z=(z_1,...,z_r).$
We keep the notation in \ref{INN}, 
\ref{IPoint} and \ref{INP}.
We note that $a_j^{(n)}\in {\mathbb D}_+\setminus \{0\},$ $\lim_{n\to\infty}a_j^{(n)}=a_j>0$
uniformly on $S,$ 
$\sum_{j=1}^n a_j^{(n)}) x_{ij}\to z_i$ 
(or $(r(j)x_{ij})v_n\to z$) uniformly on $S$ for $i=1,...,r$ as $n\to\infty$
($v_n=(\rho_A(g_1),\rho_A(g_2),...,\rho_A(g_{s(n)}))$), by \ref{ILarr}.
So Lemma \ref{ILKEY2} can be applied.
Fix $\dt>0$ and integer $K>0$ given by Lemma \ref{ILKEY2}.
We also note that since $H$ was given before we construct $q_n$
so  summands of $H$ appears in $\{\psi_j\}.$
Therefore, if $g\in {\rm ker}{\tilde \rho},$ then 
$g\in {\rm ker}\tau$ and $g\in {\rm ker} ([\Phi_1^{(0)}]\oplus [\Phi_1^{(0)}]\oplus\cdots \oplus [\Phi_1^{(0)}])$
(for any finitely many copies of $\Phi_1^{(0)}$).

Step (III). (Define $\Psi_1$ )

Let  ${\tilde \Phi}_1^{(0)}$ be a direct sum of $MK^3(k_0+1)!-1$ copies of $\Phi_1^{(0)}.$  
Note that ${\tilde \Phi}_1^{(0)}$ is a \hm\, on $A_1.$ Set $\Psi_1=\Phi_1\oplus {\tilde \Phi}_1^{(0)}.$  
If $F$ is a finite dimensional \CA\, then one has the   following
commutative diagram:
\vspace{-0.6in}
$$
\put(-85,-40){$K_0(F)$}
\put(0,-40){$K_0(F,{\mathbb Z}/k{\mathbb Z})$}
\put(105,-40){$K_1(F)$}
\put(-85, -70){$K_0(F)$}
\put(0,-70){$K_1(F, {\mathbb Z}/k{\mathbb Z})$}
\put(105,-70){$K_1(F)$}
\put(-45,-38){\vector(1,0){35}}
\put(70,-38){\vector(1,0){25}}
\put(-75,-60){\vector(0,1){15}}
\put(115, -45){\vector(0,-1){15}}
\put(-7,-68){\vector(-1,0){35}}
\put(95,-68){\vector(-1,0){25}}
\put(15, -15){}
\put(15, -92){}
\put(160, -12){}
\put(170, -88){}
$$

\vspace{-0.3in}
\noindent
where $K_0(F,{\mathbb Z}/k{\mathbb Z})=K_0(F)/kK_0(F),$
$K_1(F)=0,$ $ K_1(F, {\mathbb Z}/k{\mathbb Z})=0.$
Since $\Phi_1^{(0)}$ factors through a finite dimensional \SCA, 
it is easy to check that
$$
[\Phi_1^{(0)}]|_{K_1(A)\cap G }=0,\,\,\,\andeqn [\Phi_1^{(0)}]|_{K_1(A, {\mathbb Z}/k{\mathbb Z})\cap G}=0.
$$
Moreover,
$$
(k_0)![\Phi_1^{(0)}]|_{K_0(A, {\mathbb Z}/k{\mathbb Z})\cap G}=0\,\,\,(k\le k_0).
$$

Therefore 
$$
[\Psi_1]|_{K_1(A)\cap G}=\alpha|_{K_1(A)\cap G },\,\,\,
[\Psi_1]|_{K_1(A, {\mathbb Z}/k{\mathbb Z})\cap G}=
\alpha|_{K_1(A, {\mathbb Z}/k{\mathbb Z})\cap G },\,\,\,
$$
and $[\Psi_1]|_{K_0(A, {\mathbb Z}/k{\mathbb Z})\cap G}=
\alpha|_{K_0(A, {\mathbb Z}/k{\mathbb Z})\cap G}.$
Without loss of generality, we may 
assume that $\Psi_1(1_A)$ is a projection in $M_R(B).$

Step (IV) (Cut $\Psi_1$ according to $M,$ $K,$ $k_0$ and $\dt$)

We may assume that there are projections ${\bar p_1},...,{\bar p}_l$ in $M_k(B)$ (for some $k>0$)
such that
$$
\|{\bar p}_i-\Psi_1(p_i)\|<\eta<1/4,\,\,\,\,\,i=1,...,l.
$$
So $[{\bar p}_i]=[\Psi_1(p_i)],$ $i=1,...,l.$ 
There are projections ${\bar q_i}\le {\bar p}_i$ such that
$[{\bar q}_i]=K(k_0+1)! [\Phi_1^{(0)}(p_i)].$ Set ${\bar e_i}={\bar p}_i-{\bar q}_i.$
Note that $[{\bar e}_i]=\alpha([p_i]),$ $i=1,2,...,l.$ 
Set 
${\cal P}_1=\{{\cal P}\cup [{\tilde \Psi}_1({\cal P})]\}\cup {\cal P}_0'\cup\alpha({\cal P}),$ 
where \\
${\cal P}_0'=\{[{\bar q}_i], {\bar e_i}, [{\bar p}_i], [\Phi_1^{(0)}([p_i]):i=1,...,l\}.$
Set $G_1=G({\cal P}_1).$ We write $G_0=G_{00}\oplus G_{01},$
where $G_{00}=G_0\cap {\rm ker}\tau$ and $G_{01}\cong \tau(G_0).$ 
Let $d_1,...,d_t$ be  positive and generate $G_{01}.$

By applying \ref{ID}, 
for any $\ep>0,$ any finite subset 
${\cal G}\subset M_R(B),$ 
any $0<r<\dt<1,$  
there is a ${\cal G}$-$\ep$-multiplicative 
map $L: M_R(B)\to M_R(B)$ satisfying the following:

(1) $[L]|_{{\cal P}_1}$ and $[L]|_{G_1}$ are  well-defined and $[L]|_{G_1}$
is positive on $G_1,$ 

(2) $[L]|_{G_1\cap {\rm ker}\rho_B}={\rm id}|_{G_1\cap {\rm ker}\rho_B},$
$[L]|_{G_1\cap K_0(B, {\mathbb Z}/k{\mathbb Z})}={\rm id}|
_{G_1\cap K_0(B, {\mathbb Z}/k{\mathbb Z})},$
$[L]|_{G_1\cap K_1(B)}={\rm id}|_{G\cap K_1(B)}$ and\\ 
$[L]|_{G_1\cap K_1(B, {\mathbb Z}/k{\mathbb Z})}=
{\rm id}|_{G_1\cap K_1(B, {\mathbb Z}/k{\mathbb Z})}$
for those $k$ with 
$
G_1\cap K_i(B, {\mathbb Z}/k{\mathbb Z})\not=\emptyset 
$ {\rm (}$i=0,1${\rm )},

(3) $\rho_B\circ [L](g)\le r\rho_B(g)$ for all $g\in G\cap K_0(B),$

(4) There are $f_1,...,f_l\in K_0(A)_+$ such that
$$
\alpha(d_i)-[L](\alpha(d_i))=MK^3(k_0+1)!f_i, \,\,\,i=1,2,...,t.
$$
(Note that we may identify $\alpha$ with identity map.)
We choose $r$ so that 
$$
R\cdot r\rho_B([\Phi_1]([p_i]))\le (\dt/2MK^3(k_0+1)!) \rho(\alpha([p_i])).
$$ 
This is possible since 
$\tau(\alpha([p_i])=\tau([p_i])>0$ for all $i$ ($A$ and $B$ are simple)
and $T(A)$ is compact.
Therefore 
$$
\tau\circ [L]\circ[\Psi_1]([p_i])\le (\dt/2) \tau(\alpha([p_i]))
\andeqn
\alpha([p_i])-[L\circ \Psi_1]([p_i])>0
$$
for all $\tau\in T(A),$ 
$i=1,...,l.$

Let $[p_i]=\sum_j^t m_j^{(i)}d_j+s,$ where $m_j^{(i)}\in {\mathbb Z}$ and 
$s\in G_{00}.$
Then by (2) above, $\alpha(s)-[L]\circ \alpha(s)=0.$ 
Therefore 
\begin{eqnarray*}
\alpha([p_i])-[L\circ \Psi_1])([p_i])
&=&\alpha(\sum_jm_j^{(i)}d_j)-[L\circ \alpha](\sum_jm_j^{(i)}d_j)-
MK^3(k_0+1)![L\circ \Phi_1^{(0)}]([p_i])\\
&=&MK^3(k_0+1)!(\sum_jm_j^{(i)}f_j-[L]\circ [\Phi_1^{(0)}]([p_i]))=MK^3(k_0+1)!f_j'
\end{eqnarray*}
for some $f_j'\in K_0(B),$ $j=1,2,...,l.$


Since $K_0(B)$ is weakly unperforated, 
$f_i'>0$ for $i=1,2,...,l.$
Set $\beta: G({\cal P})\cap K_0(A)\to K_0(B)$ 
by 
$\beta([p_i])=K^3(k_0+1)!f_i',$ $i=1,2,...,l.$

 Step (V): (Construct $h$)
Let ${\tilde z_i}'=\beta([p_i])$ and
${\tilde z'}=({\tilde z_1}', {\tilde z_2}',..., {\tilde z_r}')$ (as a column).
Let ${\tilde z''}=((k_0+1)!f_1',...,(k_0+1)!f_r').$ Then $K^3{\tilde z''}={\tilde z}'.$ 
Set $z'_i=\rho({\tilde z_i}'),$ $z'=(z_1',...,z_r'),$ 
Then 
$$
\|Mz'-z\|_{\infty}<\dt\,\,\,\,\,\,\,\,\,\,\,\,{\rm (}z=(z_1,...,z_r) {\rm )}
$$
It follows from \ref{ILKEY2} that 
there is ${\tilde u}=(u_1,...,u_{s(n)})\in (K_0(B)^{s(n)})_+$ (for some $s(n)\ge r$) such that
$$
(r(j)x_{ij}){\tilde u}={\tilde z'}.
$$
Let $D= \psi_1(A)\oplus \psi_2(A)\oplus\cdots \oplus \psi_{s(n)}(A).$ 
Define a \hm\, $h_0: D\to M_k(B)$ (for some $k>0$)
such that $[h_0](e_j)=u_j,$ $j=1,...,s(n),$ where $e_j$ is a mininal projection 
in $\psi_j(A).$ Such $h_0$ exists since $B$ has stable rank one and 
real rank zero.
Define  $\pi: l^{\infty}({\mathbb Q})\to l^{\infty}_{s(n)}({\mathbb Q})$
by
$$
\{x_k\}\mapsto (x_1,...,x_{s(n)}).
$$
So 
$$
[h_0](\pi\circ {\tilde \rho}([p_i]))=\beta([p_i]), \,\,\,\,i=1,...,r.
$$
Since ${\rm ker}{\tilde \rho}\subset {\rm ker}\tau\circ \alpha\cap {\rm ker} H,$
if $x\in {\rm ker}{\tilde \rho}, $ then  
$[L\circ \Psi_1](x)=[L]\circ \alpha(x).$ 
Recall that we identify $K_0(A)$ with $K_0(B).$
By (2) above, since $\alpha(x)\in {\rm ker}\rho_B,$
$[L]\circ \alpha(x)=\alpha(x).$
Hence 
$$
\alpha(x)-[L\circ \Psi_1](x)=0.
$$
Therefore we may view that $\alpha-[L\circ \Psi_1]$ gives a \hm\,
on ${\tilde \rho}(G_0).$ 
%
%
Also $Mg$ is in the subgroup generated 
by ${\tilde \rho}([p_1]),...,{\tilde \rho}([p_r])$ for any $g\in {\tilde \rho}(G_0).$
Combining these two facts, we obtain
$$
M[h_0](\pi\circ {\tilde \rho}([p_i]))=M\beta([p_i]),\,\,\, i=1,...,r,...,l.
$$
Set $h_0'=h_0\oplus\cdots h_0,$ $M$ copies of $h_0.$
Then 
$$
[h_0'](\pi\circ {\tilde \rho}([p_i]))=\alpha([p_i])-[L\circ \Psi_1]([p_i]),\,\,\, i=1,...,r,...,l.
$$
Set $h=h_0'\circ (\psi_1\oplus \psi_2\oplus\cdots\oplus \psi_{s(n)}).$ Then $h$ is 
also ${\cal F}$-$\eta$-multiplicative.

Moreover,
$$
[h_0']\circ [\psi_1\oplus \psi_2\oplus\cdots\oplus \psi_{s(n)}]([p_i])
=[h_0']\circ (\pi\circ {\tilde \rho}([p_i]))\,\,\,\,i=1,2,...,l.
$$ 
We also have 
$$
[h]|_{{\cal P}\cap K_1(A)}=0\andeqn [h]|_{{\cal P}_0}=\alpha|_{{\cal P}_0}-[L]\circ [\Psi_1]|)|_{{\cal P}_0}.
$$
Note that, since $[h]([p_i])=MK^3(k_0+1)! f_i',$ $i=1,...,r.$
$[h](g)=MK^3(k_0+1)!g'$ for any $g\in G_0$ and some $g'\in K_0(B).$ 
Note since $D$ is finite dimensional,
we have the following (exact) commutative diagram.
\vspace{-0.6in}
$$
\put(-85,-40){$K_0(D)$}
\put(0,-40){$K_0(D,{\mathbb Z}/k{\mathbb Z})$}
\put(105,-40){$K_1(D)$}
\put(-85, -70){$K_0(D)$}
\put(0,-70){$K_1(D, {\mathbb Z}/k{\mathbb Z})$}
\put(105,-70){$K_1(D)$}
\put(-45,-38){\vector(1,0){35}}
\put(70,-38){\vector(1,0){25}}
\put(-75,-60){\vector(0,1){15}}
\put(115, -45){\vector(0,-1){15}}
\put(-7,-68){\vector(-1,0){35}}
\put(95,-68){\vector(-1,0){25}}
\put(15, -15){}
\put(15, -92){}
\put(160, -12){}
\put(170, -88){}
$$
Again, $K_1(D)=0,$  
$K_1(D, {\mathbb Z}/k{\mathbb Z})=\{0\}$ and
$K_0(D, {\mathbb Z}/k{\mathbb Z})=K_0(D)/kK_0(D)$ 
for all $k.$ 
Since $h$ factors through $D,$ 
and $[h](g)=(k_0+1)!M K^3g'',$
we conclude that
$$
[h]|_{G\cap K_1(A, {\mathbb Z}/k{\mathbb Z})}=0\andeqn
[h]|_{G\cap K_0(A, {\mathbb Z}/k{\mathbb Z})}=0\,\,\,\,{\rm for}\,\,\, k\le k_0.
$$
This implies that 
$$
[h]|_{\cal P}=\alpha|_{\cal P}-[L]\circ [\Psi_1]|_{\cal P}.
$$

Now define 
$H_1=L\circ \Psi_1\oplus h.$
Then $H_1$ is ${\cal F}$-$\eta$-multiplicative and 
$$
[H_1]|_{\cal P}=[h]|_{\cal P}+[L]\circ [\Psi_1]|_{\cal P}=\alpha|_{\cal P}.
$$
Finally, since $[H_1(1_A)]=[1_B],$
by conjegating a unitray in $(B\otimes {\cal K}{\tilde )},$
we may assume that $H_1(1_A)=1_B.$

\QED

\section{The Main Theorem}

We will use the following uniqueness theorem:

\begin{thm}{\rm (Theorem 2.3 in \cite{Lncltaf})}\label{IIIL3}
Let $A$ be a separable unital nuclear simple
\CA\, with $TR(A)=0$ satisfying the UCT.
Then, for any $\ep>0,$ and any finite subset ${\cal F}\subset A,$
there exist $\dt>0,$ a finite subset ${\cal P}\subset
P(A)$ and a finite subset ${\cal G}\subset A$ satisfying the following:
for any unital \CA \, $B$ of real rank zero and stable rank one
with weakly unperforated $K_0(B),$ and
any two ${\cal G}$-$\dt$-multiplicative morphisms $L_1, L_2: A\to B$
with
$$
[L_1]|_{\cal P}=[L_2]|_{\cal P},
$$
there exists a unitary $U\in B$ such that
$$
ad(U)\circ L_1\approx_{\ep} L_2
\,\,\,{\rm on}\,\,\,{\cal F}.
$$
\end{thm}

\begin{thm}\label{IITM}
Let $A$ and $B$ be two unital separable simple
nuclear \CA s with 
$TR(A)=TR(B)=0$ which satisfy 
the UCT.
Then 
$A\cong B$ if and only if 
$$
(K_0(A), K_0(A)_+, [1_A], K_1(A))
\cong (K_0(B), K_0(B)_+, [1_B], K_1(B)).
$$
\end{thm}

{\it Proof}: The ``only if" part is obvious and known.
We need to prove the ``if " part only.
Using the terminology in \cite{Lncltaf},
Theorem \ref{IILEX} implies that both $A$ and $B$ are pre-classifiable.
Then, by \ref{IIIL3}, Theorem 3.7 in \cite{Lncltaf} says 
that $A\cong B$ if 
$$
(K_0(A), K_0(A)_+, [1_A], K_1(A))
\cong (K_0(B), K_0(B)_+, [1_B], K_1(B)).
$$
\QED


{\small
\begin{thebibliography}{BKR}

\bibitem[Bl]{Bl}  B. Blackadar, {\em $K$-Theory for Operator Algebras},
Springer-Verlag, New York, (1986).

\bibitem[BH]{BH} B. Blackadar and D. Handelman, {\em Dimension functions and
traces on \CA s}, J. Funct. Anal., {\bf 45} (1982), 297-340.

\bibitem[BKR]{BKR}  B. Blackadar, A. Kumjian and M. R\o rdam,
{\em Approximately central matrix units and the structure of
noncommutative tori}, K-theory {\bf 6} (1992), 267-284.

 \bibitem[BBEK]{BBEK} B. Blackadar, O. Bratteli, G. A. Elliott 
 and  A. Kumjian, {\em Reduction of real rank in inductive limits of $C\sp *$-algebras},
 Math. Ann. {\bf 292} (1992),
111--126.

\bibitem[BK1]{BK1} B. Blackadar and E. Kirchberg, {\em Generalized inductive
limits of finite-dimensional \CA s},
 Math. Ann. {\bf 307} (1997), 343-380.

 \bibitem[BK2]{BK2} B. Blackadar and E. Kirchberg, {\em Inner
 quasidiagonality and strong NF algebras}, preprint.

 \bibitem[BEEK]{BEEK} O. Bratteli, G. A. Elliott, D. Evans and A. Kishimoto,
 {\em On the classification of $C\sp *$-algebras of real rank zero, III}, The infinite case.
Operator algebras and their applications, II (Waterloo, ON, 1994/1995), 11--72, 
Fields Inst. Commun., {\bf 20}, Amer. Math. Soc., Providence, RI, 1998. 


\bibitem[BP]{BP} L. G. Brown and G. K. Pedersen, {\em \CA s of
real rank zero},
J. Funct. Anal. {\bf 99}(1991), 131-149.

\bibitem[CE]{CE} M-D. Choi and E. Effros, {\em The completely positive
lifting problem for \CA s}, Ann. Math., {\bf 104} (1976), 585-609.

\bibitem[D1]{D1} M. Dadarlat, {\em Reduction to dimension three
of local spectra of real rank zero \CA s}, J. Reine Angew. Math.,
{\bf 460} (1995), 189-212.


 \bibitem[D2]{D3} M. Dadarlat, {\em Residually finite-dimensional 
 $C\sp *$-algebras}, Operator algebras and operator theory (Shanghai, 1997), 
 45--50, Contemp. Math., 228,
Amer. Math. Soc., Providence, RI, 1998.

\bibitem[DNNP]{DNNP} M. Dadarlat, G. Nagy, A. Nemethi and C. Pasnicu,
{\em Reduction of topological stable rank in inductive limits of $C\sp *$-algebras,}
Pacific J. Math. {\bf 153} (1992), 267--276.

\bibitem[DE1]{DE1} M. Dadarlat and S. Eilers,
{\em Approximate homogeneity is not a local property}, preprint 1997.


\bibitem[DE2]{DE2} M. Dadarlat and S. Eilers 
{\em On the classification of nuclear \CA s}, preprint May 1999.


\bibitem[DG]{DG} M. Dadarlat and G. Gong, {\em 
A classification result for approximately homogeneous 
\CA s of real rank zero}, Geom. Funct. Anal. {\bf 7} (1997), 646-711.


\bibitem[DL1]{DL1} M. Dadarlat and T. Loring,
{\em A universal multi-coefficient theorem for the Kasparov groups},
Duke J. Math., {\bf 84} (1996), 355-377.
 
\bibitem[DL2]{DL2} M. Dadarlat and T. Loring,
{\em Classifying \CA s via ordered, mod-p $K$-theory}, 
Math. Ann. {\bf 305} (1996), 601-616.

\bibitem[Ef]{Ef} E. Effros, {\em Dimensions and \CA s}, CBMS Regional 
Conf. Ser. in Math. no. 46, Amer. Math. Soc., Providence, R. I., 1980.

\bibitem[Ell1]{Ell1} G. A. Elliott, {\em On the classification of \CA s
of real rank zero}, J. Reine Angew. Math. {\bf 443} (1993), 179-219.

\bibitem[Ell2]{Ell2} G. A. Elliott, {\em The classification problem for
amenable \CA s}, Proc. ICM '94, Zurich, Switzerland, Birkhauser
Verlag, Basel, Switzerland, 922-932.

 \bibitem[Ell3]{Ell3} G. A. Elliott, 
{\em A classification of certain
 simple $C\sp *$-algebras}. Quantum and non-commutative analysis
(Kyoto, 1992), 373--385, Math. Phys. Stud.,
Kluwer Acad. Publ., Dordrecht, 1993.


\bibitem[Ell4]{Ell4} G. A. Elliott, {\em A classification of certain 
 simple $C\sp *$-algebras, II}, J. Ramanujan Math. Soc. {\bf 12} (1997),
 97--134.


\bibitem[EE]{EE} G. A. Elliott and D. E. Evans, {\em The structure of
irrational rotation \CA s}, Ann. Math. {\bf 138} (1993), 477-501.

\bibitem[EG1]{EG1} G. A. Elliott and G. Gong, {\em 
On inductive limits of matrix algebras over the two-torus}, Amer. J. Math. 
{\bf 118} (1996), 263--290.

\bibitem[EG2]{EG2} G. A. Elliott and G. Gong,
{\em On the classification of \CA s of real rank zero, II}, Ann. Math.
{\bf 144} (1996), 497-610.

\bibitem[EGL]{EGL} G. A. Elliott, G. Gong and L. Li, 
{\em On the classification of simple inductive limit \CA s, II: 
The isomorphism theorem}, preprint.

\bibitem[EGLP]{EGLP} G. A. Elliott, G. Gong,
H. Lin and C. Pasnicu, {\em Abelian $C^*$-subalgebras of \CA s of real
rank zero and inductive limit \CA s}, Duke Math. J., {\bf 85}, (1996),
511-554.

 \bibitem[EGS]{EGS}  G. A. Elliott, G. Gong and H. Su, {\em On the classification 
 of $C\sp *$-algebras of real rank zero, IV}. Reduction to local spectrum of 
 dimension
two. Operator algebras and their applications, II (Waterloo, ON, 1994/1995), 
73--95, Fields Inst. Commun., 20, Amer. Math. Soc., Providence, RI, 1998. 



\bibitem[ER]{ER} G. A. Elliott and M. R\o rdam
{\em Classification of certain infinite simple $C\sp *$-algebras. II.}
Comment. Math. Helv. 
{\bf 70} (1995),  615--638.

 \bibitem[ES]{ES} G. A. Elliott and H. Su, {\em $K$-theoretic classification for 
 inductive limit ${\mathbb Z}/2{\mathbb Z}$ actions on AF algebras}, Canad. J. Math., 
 {\bf 48} (1996), 946--958.

\bibitem[G1]{G1} G. Gong, {\em 
On the inductive limits of matrix algebras over higher
dimensional spaces}, Part I \& II: Math. Scand. {\bf 80} (1997)
40-55 \& 56-100.

\bibitem[G2]{G2} G. Gong, {\em
On the classification of simple inductive limit \CA s, I; The reduction
theorem}, preprint.

\bibitem[Go]{Go} K. R. Goodearl, {\em Notes on a class of simple $C^*$-algebras with 
real rank zero}, Publ. Mat.
(Barcelona) 36 (1992), 637-654.

\bibitem[K]{K} E. Kirchberg {\em 
The classification of purely infinite \CA s using Kasparov's theory},
Fields Institute Comm., to appear.

 \bibitem[Ki]{Ki} A. Kishimoto, {\em Non-commutative shifts and crossed products},
 preprint June 2000.

 \bibitem[JS]{JS} X. Jiang and H. Su, {\em On a simple unital projectionless 
 $C\sp *$-algebra}, Amer. J. Math., {\bf 121} (1999), 359--413


\bibitem[Li1]{Li1} L. Li {\em 
Classification of simple $C\sp *$-algebras: inductive limits of matrix algebras 
over trees}, Mem. Amer. Math. Soc. {\bf 127} (1997), no. 605.

\bibitem[Li2]{Li2} L. Li, {\em Simple inductive limit $C\sp *$-algebras: 
spectra and approximations by interval algebras},  J. Reine Angew. Math. 
{\bf 507} (1999), 57--79.


\bibitem[Ln1]{Ln2} H. Lin, {\em 
Approximation by normal elements with finite spectra in 
$C\sp *$-algebras of real rank zero}i,
 Pacific J. Math. {\bf 173} (1996),443--489.

 \bibitem[Ln2]{Ln3} H. Lin, {\em Almost commuting unitaries and classification of 
 purely infinite simple 
$C\sp *$-algebras},  J. Funct. Anal. {\bf 155} (1998), 1--24.

 \bibitem[Ln3]{Ln4} H. Lin, {\em 
On the classification of $C\sp *$-algebras of real rank zero with zero 
$K\sb 1$} J. Operator Theory {\bf 35} (1996),  147--178.

 \bibitem[Ln4]{Ln5} H. Lin, {\em
 A classification theorem for infinite Toeplitz algebras },
 Operator algebras and operator theory (Shanghai, 1997), 219--275, Contemp. Math., 228,
Amer. Math. Soc., Providence, RI, 1998. 

\bibitem[Ln5]{Ln6} H. Lin {\em 
Classification of simple \CA s with unique traces},
Amer. J. Math., {\bf 119} (1997), 1263-1289.

\bibitem[Ln6]{Lnsu} H. Lin, {\em Stable approximate unitarily equivalence
of \hm s}, J. Operator Theory, to appear..

\bibitem[Ln7]{Lntaf} H. Lin, {\em Tracially AF \CA s},
Trans. Amer. Math. Soc., to appear.

\bibitem[Ln8]{Lncltaf} H. Lin, {\em Classification of simple 
TAF \CA s}, Canad. J. Math., to appear.

\bibitem[Ln9]{LnI} H. Lin, {\em Locally type I simple tracially AF \CA s},
preprint 98.

\bibitem[Ln10]{Lntr} H. Lin, {\em Tracial topological rank 
of \CA s}, Proc. London Math. Soc., to appear.

\bibitem[Ln11]{LnQ}  H. Lin, {\em Quantized deformation of topological spaces ---Some
applications of classification of nuclear \CA s}, preprint 99.

\bibitem[Ln12]{Lntor} H. Lin, {\em 
Classification of simple \CA s and higher dimensional non-commutative tori}, 
preprint Feb. 2000.

 \bibitem[LP1]{LP1} H. Lin and N. C. Phillips, {\em Classification of direct limits of
even Cuntz-circle algebras}, Memoirs Amer. Math. Soc.,  {\bf 118}, no. 565,
(1995).


 \bibitem[LP2]{LP2} H. Lin and N. C. Phillips, {\em Approximate Unitary equivalence
of homomorphisms from $O_{\infty}$},  J. Reine Angew. Math.,
{\bf 464} (1995), 173-186.

 \bibitem[LP3]{LP3} H. Lin and N. C. Phillips, {\em Almost multiplicative morphisms and
 Cuntz algebra $O_2$},  Inter. J. Math., 6 (1995) 625-643.


\bibitem[LS]{LS} H. Lin and H. Su, {\em Classification of direct limits of 
generalized Toeplitz algebras}, Pacific J. Math. 181 (1997), no. 1, 89--140.

\bibitem[LqP]{LqP} Q. Lin and N. C. Phillips, {\em 
\CA s of minimal diffeomorphisms}, preprint 2000.


\bibitem[NT]{NT} K-E. Nielsen and Thomsen,
{\em Limits of circle algebras}, Exposition. Math. {\bf 14} (1996), 17--56.

\bibitem[Pa]{Pa} V. I. Paulsen, {\em Completely bounded maps and dilations},
Pitnam Research Notes in Math., {\bf 146}, Sci. Tech. Harlow, (1986).



\bibitem[Pd]{Pd} G. K. Pedersen, {\em \CA s and their Automorphism Groups},
Academic Press, London/New York/San Francisco, 1979.

\bibitem[Po]{Po} S. Popa, {\em On locally finite dimensional approximation of
\CA s}, Pacific J. Math., (1997), 141-158.

\bibitem[P1]{P1} N. C. Phillips {\em A classification theorem
for nuclear purely infinite simple \CA s},
Doc. Math. 5 (2000), 49--114. 

 \bibitem[P2]{P2} N. C. Phillips {\em The real rank of 
 direct limits of recursive subhomogeneous \CA s},
  in preparation.

\bibitem[Ro1]{Ro1} M. R\o rdam, {\em Classification of inductive 
limits of even Cuntz algebras}, J. Reine Abgew. Math., {\bf 440}
(1993), 175-200.

\bibitem[Ro2]{Ro2} M. R\o rdam, {\em Classification of
Cuntz-Krieger algebras}, $K$-theory,



\bibitem[Ro3]{Ro3} M. R\o rdam, {\em Classification of certain
   infinite simple \CA s}, J. Funct. Anal., {\bf 131} (1995), 415-458.


\bibitem[RS]{RS} J. Rosenberg and C. Schochet, 
{\em The Kunneth theorem and the universal Coefficient theorem for
Kasparov's generalized functor}, Duke Math. J. {\bf 55} (1987),
431-474.

\bibitem[Sc1]{Sc1} C. Schochet, {\em Topological methods for
\CA s III: axiomatic homology}, Pacific J. Math. 
{\bf 114} (1984), 399-445.

\bibitem[Sc2]{Sc2} C. Schochet, {\em Topological methods for
\CA s IV: mod p homology}, Pacific J. Math.,
{\bf 114} (1984), 447-468.


 \bibitem[Su1]{Su1} H. Su, {\em On the classification of $C\sp *$-algebras 
 of real rank zero: inductive limits of matrix algebras over non-Hausdorff 
 graphs}, Mem. Amer. Math.
Soc. {\bf 114} (1995), no. 547, 


\bibitem[Su2]{Su2} H. Su, {\em Classification for certain simple real rank zero 
$C\sp *$-algebras}. J. Funct. Anal. 140 (1996), no. 2, 505--540

\bibitem[Th1]{Th1} K. Thomsen, {\em On isomorphisms of inductive
limits of \CA s}, Proc. Amer. Math. Soc. {\bf 113} (1991), 947-953.

 \bibitem[Th2]{Th2}  K. Thomsen, {\em Inductive limits of interval algebras: 
 the tracial state space}, Amer. J. Math. 116 (1994), no. 3, 605--620

\bibitem[Zh1]{Zh1} S. Zhang, {\em \CA s with real rank zero
and the internal structure of their corona and multiplier
algebras, Part I}, Pacific J. Math., {\bf 155}(1992), 169-197.

\bibitem[Zh2]{Zh2} S. Zhang, {\em A property of purely infinite simple
\CA s}, Proc. Amer. Math. Soc., {\bf 109} (1990), 717-720.

\bibitem[Zh3]{Zh3} S. Zhang, {\em Matricial structure and homotopy type
of simple \CA s with real rank zero}, J. Operator Theory, {\bf 26} (1991),
283-312.

\end{thebibliography}
}
\end{document}